\documentclass{article}		
	\usepackage{graphicx}
	\pdfoutput=1

	\usepackage[cp1251]{inputenc}
	\usepackage[russian]{babel}
	\usepackage{amsfonts, amsmath, amsthm}
	
	\newtheorem{Th}{Теорема}
	\newtheorem{Lemm}{Лемма}
	\newtheorem{Foll}[Lemm]{Следствие}
	\newtheorem{sSt}{Предложение}[section]
	\newtheorem{ssSt}{Предложение}[subsection]
		\theoremstyle{definition}
	\newtheorem{Meth}{Метод}
	\newtheorem{sDef}{Определение}[section]
	\newtheorem{ssDef}{Определение}[subsection]
	
\newenvironment{Proof} 		
	{		
		\parskip=-3pt
		\parindent=12pt
		\par
		{\bf Доказательство.}
	}	
	{		
		\hfill
		$\square$
	}	
	
\newenvironment{Des}[1]		
	{
		{\bf Обозначения:}
		\par
		\vspace{5pt}
		\hbox{
			\hspace{30pt}
			\parbox{0.93\textwidth}%
			{#1}
		}
	}
	{
		\par
	}
	
\newenvironment{Indention}[1]	
	{\hbox{%
	\hspace{30pt}\parbox{0.93\textwidth}%
	{#1}}}{\par}
	
\begin{document}			

	
	\begin{center} 		
		\large\textbf{О числе запретов, задающих периодическую последовательность.} 
	\end{center}

\section{Предисловие}

	При изучении алгебраических объектов используется их представление в
виде образующих и определяющих соотношений. Способов представлять
элемент через образующие может быть много, так что изучается {\it
каноническая форма} представления. Например, пусть $A$~--
ассоциативная алгебра, $a_1,\dots,a_s$~-- ее образующие.  Порядок
$a_1\prec\cdots\prec a_s$ индуцирует порядок на множестве  мономов от
$a_i$ (сперва по длине, потом лексикографически). Множество слов, не
являющихся линейной комбинацией меньших, образует {\it нормальный
базис} алгебры $A$. {\it Функция роста} $V_A(n)$ есть размерность
пространства, порожденного мономами степени не выше $n$ и совпадает с
числом элементов нормального базиса степени не выше $n$. Функции роста
и нормальные формы определяются для разных алгебраических систем, им
посвящена обширная литература. Обзор -- см. \cite{BBL}

	Рассматривая идеал соотношений $I$, рассматривают старшие члены его
элементов или {\it редуцируемые слова}. Пусть $\{f_i\}$~-- образующие $I$.
Тогда надслово редуцируемого слова редуцируемо. {\it Обструкцией}
называется минимальное редуцируемое слово, т.е. без редуцируемых
подслов. Если все обструкции входят в старшие члены базиса, то он
называется {\it базисом Гребнера-Ширшова} идеала $I$. Впервые это
понятие было введено А.И.Ширшовым. Аналогичные понятия для
полилинейных слов ввел В.Н.Латышев (минимальность понимается также и в
том что слово не является изотонным образом редуцированного).
А.И.Ширшов ввел понятие {\it композиции} и предложил критерий того,
что $\{f_i\}$ является базисом гребнера. Это легло в основе знаменитой
Diamond-леммы Бергмана.

	Базис гребнера даже конечно определенной алгебры может быть
бесконечен. А.Я.Белов ввел понятие {\it короста} в алгебе или функции,
$B(n)$ выражающей количество обструкций длины не выше $n$.

Данная работа посвящена исследованию {\it кодлины} периода или
количеству запретов, которыми можно задать периодическую
последовательность.

Основной результат данной работы состоит в следующем: 

\textbf{Теорема 1.} \textit{В случае двухсимвольного алфавита $A = \{a, b\}$ если слово из $I_n$ задается $c$ запретами, то $\varphi_c\geq n$, где $\varphi_c$ --- это $c$-е число Фибоначчи ($\varphi_1=1$, $\varphi_2=2$, $\varphi_3=3$, $\varphi_4=5$ и т.д.).}

Отметим, что логарифмическая оценка, а также основные примеры и оценка сверху на количество запретов были получены ранее в работе Г. Р. Челнокова \cite{Chelnokov}. Мы вычисляем более точную ассимптотику (множитель при логарифме).

В основе доказательства лежит работа с {\it графами и схемами Рози}

Последовательности схем Рози исследовались в ряде работ. С их помощью
удается решить ряд алгоритмических проблем. см. например \cite{Mitrofanov1}, \cite{Mitrofanov2} и \cite{Mitrofanov3}

Отметим,что результаты работ \cite{Mitrofanov1} и \cite{Mitrofanov2} другим методом независимо получены
Ф.Дюрандом: \cite{Durand1}, \cite{Durand2}

 В терминах размеченных схем Рози удается получить критерий того, что
слово отвечает перекладыванию отрезков \cite{BelovChernInt}. Подробнее см. обзор \cite{BelovKondakovMitrofanov}

Перейдем к доказательству основного результата.	

\section{Введение}		

	В данной работе изучаются свойства бесконечных в обе стороны периодических последовательностей (слов, строк --- будем считать эти термины эквивалентными) над фиксированным конечным алфавитом $A$. 
	
	\begin{Des}{ 		
		$F$ --- множество всех конечных слов в алфавите $A$.\\
		$F_n$ --- слова длины $n$.\\
		$U_n$ --- слова длины не больше $n$.\\
		$F_w$ --- все конечные подслова слова $w$.\\
		$S_w =$ $F\backslash F_w$ (все конечные слова, не встречающиеся в $w$).\\
		$I$ --- бесконечные в обе стороны слова в алфавите $A$ заданные с точностью до сдвига (формально: реализованные, например, как отображения из $\mathbb Z$ в $A$ факторизованные по отношению эквивалентности сдвига).\\
		$I_n\subset I$ --- слова c наименьшим периодом n (которые можно себе представлять как n букв стоящие по кругу).\\
		$I_\infty\subset I$ --- непериодические слова. (Очевидно, $I = \left(\bigcup\limits_{n = 1}^\infty I_n\right) \cup I_\infty$)}
	\end{Des}
		
	\begin{sDef} 		
		Будем называть любое множество $S\subseteq F$ системой запретов, а его элементы --- запретами. 
	\end{sDef}
	
	\begin{sDef} 		
		Будем говорить, что бесконечное слово $w\in I$ удовлетворяет системе запретов S если $\forall s\in S\ s$ не является подсловом $w$, или, что то же самое, $S\subseteq S_w$
	\end{sDef}
	
	\begin{sDef} 		
		Будем говорить, что система запретов $S$ задает бесконечное слово $w\in I$ если оно и только оно удовлетворяет этой системе запретов.
	\end{sDef}
	
	\begin{sSt}  		
		\label{Fin->Periodic}
		Если система запретов $S$ задает какое-то слово $w\in I$ и $|S|<\infty$ то $w$ --- периодическое.
	\end{sSt}
		
	\vspace{-12pt} Это будет доказано ниже.
		
	\begin{sSt}  		
		\label{Periodic->existsFin}
		$\forall n\ \forall w \in I_n\ \exists S,\ |S|<\infty$ и $S$ задает $w$.
	\end{sSt}
		
	\vspace{-12pt} Например, $S = S_w\cap U_{n+1}$. И даже $S = S_w\cap F_{n+1}$. Это также будет доказано ниже.
		
	\vspace{10pt} В данной работе исследуется вопрос о наименьшем возможном количестве элементов в конченой системе запретов, задающей какое-нибудь слово из $I_n$ в зависимости от $n$. Доказана
	
	\begin{Th} 			
	\label{mainth}
		В случае двухсимвольного алфавита $A = \{a, b\}$ если слово из $I_n$ задается $c$ запретами, то $\varphi_c\geq n$, где $\varphi_c$ --- это $c$-е число Фибоначчи ($\varphi_1=1$, $\varphi_2=2$, $\varphi_3=3$, $\varphi_4=5$ и т.д.).
	\end{Th}
	
	\vspace{-12pt} Также рассматривается вопрос минимальности этой оценки для каждого $n$ и случай $k$-буквенного алфавита.

\section{Доказательство оценки. Теоретическая часть.}
\subsection{Определение и простейшие свойства графов Рози}

	\begin{ssDef}		
		$k$-м графом Рози слова $w\in I$ называется граф $G_k^w$ множеством вершин которого является множество $V = F_w\cap F_k$, а множеством ребер --- множество $E = F_w\cap F_{k+1}$, где ребро $e\in E$ ведет из вершины $v_1$ равной первым $k$ буквам $e$ в вершину $v_2$ равную последним $k$ буквам $e$ --- это естественная ориентация на $G_k$:
	\end{ssDef}

	\includegraphics[height=25mm]{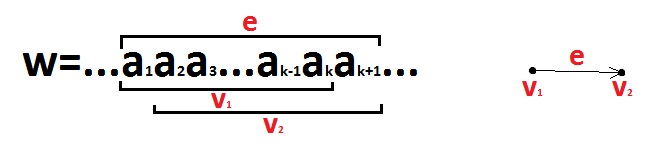}
	
	В случае, если понятно о графе Рози какого слова идет речь, верхний индекс будет опускаться. Чтобы не путаться, введем
	\begin{Des}{ 		
	$l(v)$ или $l(e)$ --- слово, соответствующее вершине $v$ или ребру $e$ соответственно,\\
	$v(l)\in G_{|l|}$ и $e(l)\in G_{|l-1|}$ --- ребро и вершина для данного слова.}
	\end{Des}
	Для некоторой конечной строки $l\in F$ 
	\begin{Des}{		
	$l[i]$ --- $i$-я буква $l$\\
	$l[i:j]$ --- строка состоящая из букв $l$ с $i$-й по $j$-ю в том же порядке.}
	\end{Des}
	Заметим, что переход от одной вершины к другой в направлении ориентации соединяющего их ребра соответствует переходу к рассмотрению подслова, сдвинутого на одну букву. Это позволяет понять что по путям в графах Рози можно выписывать слова, а по словам --- строить пути в некотороых графах Рози. В частности: 
	\begin{ssSt} 		
		Любому замкнутому пути $e_1\cdots e_m$ в некотором графе Рози $G_k^w$ соответствует слово $w_1\in I_m$: $w_1 = \cdots l(e_m)[1] + l(e_1)[1] + l(e_2)[1] + \cdots + l(e_m)[1] + l(e_1)[1]\cdots$. Поскольку ребра последовательные, то все $l(e_i)$ встречаются как подслова в $w_1$ (начинающиеся в соотвествтующем $l(e_i)[1]$): $l(e_i) = l(e_i)[1] + l(e_{i+1})[1]\cdots l(e_{(i+k) \mod m})[1]$.
	\end{ssSt}
	
	Также с точностью до сдвига выписывается бесконечное слово (возможно, непериодическое) по любому бесконечному пути в $G_k^w$. С точностью до сдвига --- т.к. кроме первой можно выписывать любую наперед заданную по порядку букву каждого ребра.
	
	\begin{ssSt} 		
		\label{Path}
		В любом графе Рози $G_k^w$ есть бесконечный в обе стороны путь проходящий по всем ребрам и всем вершинам, соответствующий слову $w$.
	\end{ssSt}
	\begin{Proof} 
	Выберем в слове $w\ k+1$ букву идущую подряд (подслово $l_1$). $e(l_1)$ --- первое ребро пути. Cместим выбор на 1 букву вправо (подслово $l_2$). $e(l_2)$ --- это второе. И так далее. Аналогично, влево. Очевидно, конец $e(l_i)$ совпадает с началом $e(l_{i+1})$. Очевидно этот путь проходит по всем ребрам и всем вершинам: он перебирает все подслова слова $w$ длин $k$ и $k+1$. 
	\end{Proof}
	
	В случае периодического слова $l_{i+n} = l_i$, и $e(l_{i+n}) = e(l_i)$, и можно рассматривать просто замкнутый путь длины $n$, т.к. остальное --- его повторение.
	
	\begin{ssDef}		
		Будем называть систему запретов $S$ приведенной системой запретов, если она задает некоторое слово $w\in I$, не содержит дубликатов и любое из собственных подслов любого слова $s\in S$ является подсловом $w$.
	\end{ssDef}
	
	\begin{ssSt} 		
		\label{Reduction}
		По любой конечной системе запретов $S$, задающей некоторое бесконечное слово $w$ можно построить приведенную приведенную систему запретов задающую то же слово и содержащую не большее число элементов. 
	\end{ssSt}
	\begin{Proof} 
	Если какое-то собственное подслово $s_1$ какого-то слова $s\in S\subseteq S_w$ так же не содержится в $w$ --- можно заменить $s$ на $s_1$ и удалить дубликаты. Очевидно, $w$ подходит новой системе (никакое слово в нем не содержится) и если какое-то слово подходит новой системе, то оно подходило и старой (т.к. запреты мы только укоротили и не добавляли) --- а значит новая система тоже задает $w$. В силу конечности количества запретов и их длин --- процесс когда-нибудь завершится и мы получим приведенную систему запретов. 
	\end{Proof}

	Изучим, как связанны $G_k^w$ и $G_{k + 1}^w$ для некоторого периодического $w\in I_n$. Понятно, что вершины $G_{k + 1}$ взаимно соответствуют ребрам $G_k$ --- т.к. это подслова $w$ длины $k + 2$. А ребра $G_{k + 1}$ соответствуют путям длины 2 в $G_k$ (парам последовательных ребер): ребру $e$ в $G_{k+1}$ соответствует пара ребер $e(\ l(e)[1:k+1]\ )$ и $e(\ l(e)[2:k+2]\ )$ --- с общей вершиной $v(\ l(e)[2:k+1]\ )$.
	\begin{ssSt} 		
		\label{obligatoryRestrict}
		А тем путям длины 2 в $G_k$, которым не соответствует ребро графа $G_{k+1}$, соответствует запрет в любой приведенной системе запретов $S$, задающей данное слово, точнее --- существует биекция между запретами из $S$ длины $k+1$ и такими путями в $G_k$, причем каждой такой паре последовательных ребер $e_1$ и $e_2$ соответствует запрет $l(e_1) + l(e_2)[k+1]$, а запрету $s$ --- пара ребер $e(\ s[1:k+1]\ )$ и $e(\ s[2:k+2]\ )$.
	\end{ssSt}
	\begin{Proof} 
	Предположим, что найдутся два последовательно идущих ребра $G_k$\ $e_1$ и $e_2$ таких, что слово $x = l(e_1) + l(e_2)[k+1]$ не является подсловом $w$ и не принадлежит $S$. Значит, никакой запрет из $S$ не содержит $x$ как подслово. Пусть слева от $l(e_1)$ в $w$ написано бесконечное влево слово $u$, а справа от $l(e_2)$ --- бесконечное вправо слово $v$. Тогда слово $w_1 = uxv$ не равное $w$ т.к. в $w$ не встречается $x$ --- будет удовлетворять системе $S$ (если какой-то запрет из $S$ содержится в $w_1$, то либо он содержит $x$, либо он содержится в $u + l(e_1)$ или $l(e_2) + v$ --- чего быть не может). Противоречие. Для любого же запрета $s\in (S\cap F_{k+2})$, $s$ не является подсловом $w$ а $s[1:k+1]$ и $s[2:k+2]$ (поскольку система приведенная) --- являются. Значит запрету $s$ соответствует пара последовательных ребер в $G_k$, которым не соответствует ребра $G_{k+1}$.
	\end{Proof}
		
	\begin{Des}{ 		
	$\tilde{v}(e)=v(l(e))$\\
	$\tilde{e}(v)=e(l(v))$\\
	$\tilde{e}(e_1,e_2) = e(l(e_1) + l(e_2)[k+1])$\\
	$\tilde{p}(e) = (e(\ l(e)[1:k+1]\ ), e(\ l(e)[2:k+2]\ )$ --- в последних двух функциях аргументом или значением в соответствующих случаях могут выступать элементы приведенной системы запретов.}
	\end{Des}
	
	Это значит, что по $G_k$ очень легко нарисовать $G_{k+1}$:
	
	\begin{Meth} 		
		\label{Trivial} 
		В центре каждого ребра поставить по вершине, соеденить ребрами все пары последовательных ребер $G_k$, стереть все старое (от $G_k$) и зачеркнуть (точнее, стереть) все запрещенные ребра (ребра, принадлежащие $S\cap F_{k+2}$). 
	\end{Meth}
	
	\begin{ssDef}		
		Будем называть произведение первой части данных действий "`полушагом"', а второй (стирания ребер, соответствующих словам из $S$) --- "`применением запретов"'. Всю совокупность этих действий (а точнее --- переход от от рассмотрения $G_k$ к $G_{k+1}$) будем называть "`шагом"'.
	\end{ssDef}
	\begin{Des}{		
		$h(G_k)$ --- граф, получающийся из $G_k$ после полушага. (вообще, аналогичные операции можно проделать с любым ориентированным графом $G$. Результат --- тоже ориентированный граф --- будем обозначать $h(G)$). За $\tilde{p}$ будем обозначать естественное вложение вершин, ребер и развилок всех типов из $G_{k+1}$ в $h(G_k)$.}
	\end{Des}

\subsection{Развилки в графах Рози. Их количество.}
\subsubsection{Определения развилок}

	Заметим, что поскольку букв в нашем алфавите всего две, входящая и исходящая степень каждой вершины не больше двух (есть всего два способа продолжить на одну букву слово соответствующее данной вершине как влево, так и вправо). А поскольку есть замкнутый путь, проходящий через все ребра и все вершины (предложение \ref{Path}) --- эти степени не меншье одного. То есть есть четыре типа вершин (указаны соотв. входящая и исходящая степени): (1,1), (2,1), (1,2) и (2,2). 
	\begin{ssDef} 		
		Назовем такие вершины соответственно дорога, входящая развилка, исходящая развилка и перекресток. 
	\end{ssDef} 
	При полушаге они преобразуются так: 
	
	\includegraphics[height=35mm]{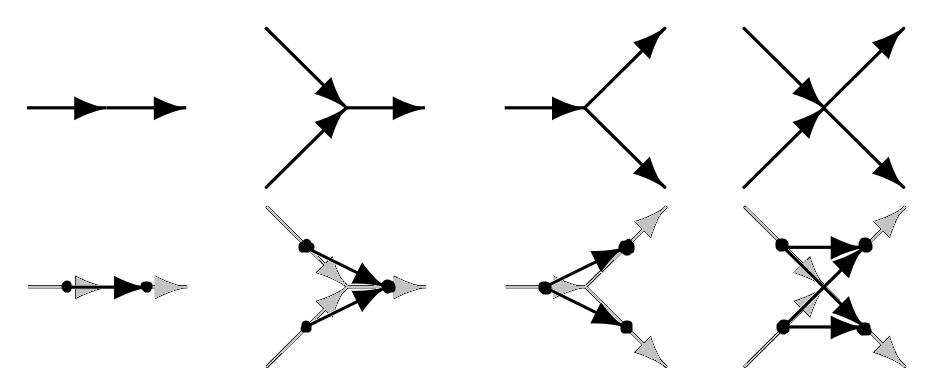} 

	В дальнейшем для нас будет важно понятие развилок. Стоит обратить внимание что можно ввести три разных (однако почти эквивалентных по свойствам) определения развилок, которые каждое по своему удобны для доказательства некоторых утверждений. Первое уже введено.
	
	\begin{ssDef}		
		Будем называть входящей $d$-развилкой пару различных ребер, имеющих общий конец, а исходящей --- общее начало. 
	\end{ssDef}
	
	\begin{ssDef}		
		Будем называть входящей $t$-развилкой пару из входящей $d$-развилки и "`направляющего"' ребра (возможно --- принадлежащего этой $d$-развилке), начало которого совпадает с ее концом , а исходящей --- наоборот.
	\end{ssDef}
	
	При использовании каждого определения будем употреблять соответственно слова "`развилка"', "`$d$-развилка"' и "`$t$-развилка"' за исключением случаев, которые будут оговорены отдельно.

	\begin{Des}{		
		$c_k^w$ --- количество перекрестков в $G_k^w$\\[2pt] 
		$f_k^w$ --- количество входящих развилок (на самом деле, входящих или исходящих --- не важно, т.к. они равны --- это будет доказано в следующем абзаце)\\[2pt] 
		$d_k^w$ --- количество входящих $d$-развилок\\[2pt] 
		$t_k^w$ --- количество входящих $t$-развилок.}
	\end{Des}
	
	Если понятно о каком слове идет речь, верхний индекс будет опускаться.
	
	\begin{ssSt}		
		$\forall w\ f_k^w$ равно количеству исходящих развилок в $G_k^w$.
	\end{ssSt}
	\begin{Proof} 
		Тривиально следует из того что сумма входящих степеней в любом $G_k^w$ равна сумме исходящих. 
	\end{Proof}

	\begin{ssSt} 		
		\label{forkNum}
		$d_k^w = f_k^w + c_k^w$, $t_k^w = f_k^w + 2*c_k^w$. Для исходящих $d$- и $t$-развилок аналогично.
	\end{ssSt}
	\begin{Proof} 
		Cопоставим каждой $d$-развилке ее конец. Его исходящая степень может быть 1 или 2. Для $t$-развилок аналогично. Для исходящих $d$- и $t$-развилок аналогично. 
	\end{Proof}

	Значит и количества входящих $d$- и $t$-развилок тоже равны соответствующим количествам исходящих.

\subsubsection{Соответствия между развилками, перекрестками, $d$- и $t$-развилками в одном и различных графах Рози}

	Как уже упоминалось, определения развилок, $d$-развилок и $t$-развилок почти эквивалентны. В частности --- между ними есть естественные соответствия: \\[4pt] 
	\begin{Indention}{	
	Развилке --- соответствуют одна $d$-развилка и одна $t$-развилка (входящей --- входящие, исходящей --- исходящие): $\tilde{d}(f)$ и  $\tilde{t}(f)$ (поскольку просто развилками мы называем вершины графов Рози, аргументами могут быть соответствующие вершины).\\[4pt] 
	Перекрестку ---  две $d$-развилки (входящая и исходящая) и четыре $t$-развилки (две входящие и две исходящие). $\tilde{d}(c)$ и  $\tilde{t}(c)$ (Аналогичное замечание. Кроме того, образом многозначного соответствия мы естественно будем считать множество из всех соответствующих элементов).\\[4pt] 
	$d$-развилке --- ее вершина (развилка соотв. типа или перекресток) и соответствуеющие ей $t$-развилки: $\tilde{v}(d)$ и  $\tilde{t}(d)$.\\[4pt] 
	$t$-развилке --- тоже: $\tilde{v}(t)$ и  $\tilde{d}(t)$.}
	\end{Indention}
	
	Заметим, что в случае развилки это соответствие взаимно однозначное, а в случае перекрестка --- нет.
	
	Заметим также, что все свойства графов Рози для $h(G_k)$ выполняются (что и понятно --- он мог бы быть графом Рози какого-нибудь слова, если бы $S$ не содержала запретов длины $k+2$). Значит можно говорить о развилках в нем. 
	
	\begin{ssSt} 		
		\label{bijectHalf}
		Существует биекция между $t$-развилками в $G_k$ и $d$-развилками в $h(G_k)$ (и значит их количества равны), причем $t$-развилке соответствует $d$-развилка, вершина которой соответствует направляющему ребру этой $t$-развилки.
	\end{ssSt}
	\begin{Proof} 
		Сопоставим естественным образом вершины и ребра $h(G_k)$ (k+1)- и (k+2)-буквенным словам (вершине поставленной в центре ребра $e$ --- слово $l(e)$, а ребру из $e_1$ в $e_2$ --- слово $l(e_1) + l(e_2)[k+1]$ . Теперь зададим биекцию: каждой $t$-развилке $t_1$ в $G_k$ состоящей из $d$-развилки $d$ (из ребер $d_1$ и $d_2$) и ребра $e_1$ будет соответствовать $d$-развилка из ребер соотв. путям $d_1e_1$ и $d_2e_1$, а $d$-развилке в $h(G_k)$ будет соответствовать $t$-развилка из пары ее начальных вершин и конечной вершины --- т.к. начальные вершины не совпадают: нигде кроме $G_0$ не может два разных ребра начинаться и заканчиваться в одинаковых вершинах --- иначе у конечной вершины последняя буква должна быть равна и "`a"', и "`b"', а такого быть не может (заметим, что с тем же успехом можно было рассмотреть первую букву начальной вершины). 
	\end{Proof}
	\begin{Des}{
		Обозначим это соответствие за $\tilde{h}(t)$ (Обратное можно обозначать и как $\tilde{h}^{-1}(d)$, и как $\tilde{h}(d)$ --- т.к. оно однозначно определяется типом аргумента)}
	\end{Des}
	
	Заметим, что если вершина $d$-развилки в $h(G_k)$ --- развилка, то ей соответствует единственная $t$-развилка. Аналогично если вершина $t$-развилки в $G_k$ --- развилка, то ей соответствует единственная $d$-развилка. Значит, на множестве всех $d$-развилок во всех $G_k$ данного слова можно ввести отношение эквивалентности порожденное следующими равенствами: две $d$-развилки эквивалентны, если они находятся в последовательных графах Рози $G_k$ и $G_{k+1}$ и однозначно друг другу соответствуют в указанном естественном смысле (если вершина той что в $G_k$ --- развилка). Аналогично для $t$-разивлок. Формально: $d_1 ~ d_2$, если  $\tilde{p}(d_2) = \tilde{h}(\tilde{t}(d_1))$, (или наоборот).
	
	\begin{ssDef} 		
		Соответствующие классы эквивалентности будем называть $md$-развилками и $mt$-развилками.
	\end{ssDef}
	В дальнейшем нам будет удобно говорить о какой-то конкретной $md$- или $mt$-развилке просто указав один из ее элементов. Заметим, что $d$- или $t$-развилки из одного класса эквивалентности содержатся в некотором количестве подряд идущих графов Рози (по 1 в каждом). Для удобства терминологии будем представлять себе $k$ как ось дискретного времени (для $G_k$) и говорить, что $md$- или $mt$-развилка "`существует в момент времени $t$"' или "`принадлежит $G_t$"', если какой-то ее элемент содержится в $G_t$.
	
	Очевидно, что существование каждой $md$- и $mt$-развилки начинается перекрестком и заканчивается перекрестком или запретом (на самом деле, заканчиваться оно может тоже только перекрестком, т.к. как будет видно далее, в случае приведенной системы запретов, запреты могут происходить только сразу через шаг после перекрестков, а значит этой развилки бы просто не существовало вообще).
	
\subsubsection{Количество развилок}

	\begin{ssSt} 		
		\label{halfDiff}
		Количество $d$-развилок при полушаге увеличивается на $c_k$.
	\end{ssSt}
	\begin{Proof} 
		Очевидно следует из предложений \ref{forkNum} и \ref{bijectHalf}. 
	\end{Proof}

	\begin{ssSt} 		
		\label{restrictDiff}
		При применении запретов количество входящих $d$-развилок уменьшается на количество примененных запретов (стертых ребер).
	\end{ssSt}
	\begin{Proof} 
		Пусть ребро $x$ с вершинами $v_1$ и $v_2$ из $G_{k+1}$ должно быть стерто (содержится в приведенной с.з. $S$). Исходящая степень $v_1$ как и входящая степень $v_2$ в $G_{k+1}$ не меньше одного. Значит в $h(G_k)$ эти степени не меньше 2. Но они не могут быть больше 2. Значит они ровно 2. Очевидно, одно ребро не может содержаться в двух разных входящих $d$-развилках. Значит количество входящих $d$-развилок при стирании этого ребра уменьшится ровно на 1. Далее индукция. 
	\end{Proof}

	\begin{Lemm} 		
		\label{forkDiff}
		$\forall w\in I\ \forall k\geq 1,\ d_k^w = d_{k-1}^w + c_{k-1}^w - |S\cap F_{n+1}|$, где $S$ --- приведенная система запретов, задающая $w$.
	\end{Lemm}
	\begin{Proof} 
		Очевидно следует из предложений \ref{halfDiff}, \ref{obligatoryRestrict} и \ref{restrictDiff}.
	\end{Proof}

	\begin{Lemm} 		
		\label{graphFin}
		Если $w\in I_n$, то $G_n^w$ --- циклический граф содержащий $n$ ребер (то есть $d_n^w = 0$ а $|G_n^w| = n$).
	\end{Lemm}
	\begin{Proof} 
		По условию $v_1 = v_2$ (см. рисунок), значит в $G_n$ есть замкнутый путь по $n$ ребрам (т.к. при $n$ сдвигах рассматриваемого подслова на одну букву мы получим то же подслово). Он не самопересекается, иначе у $w$ есть меньший период ("'$a_{2k+1}...a_{k+n}$"'$=$"'$a_{k+1}...a_n$"'$=$"'$a_1...a_{n-k}$"'; 
	
	"'$a_{k+n+1}...a_{2k+n}$"'$=$"'$a_1...a_k$"'$=$"'$a_{n-k+1}...a_n$"') --- равный $k$ (где $k$ это количество ребер до самопересечения, а $v_1$ --- вершина, в которой путь самопересекается). 
	\end{Proof}
	
	\includegraphics[height=25mm]{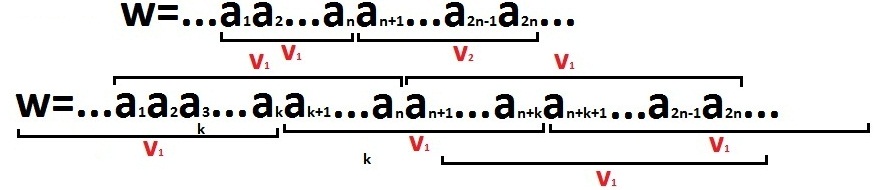}

\subsubsection{Долг.}		

	Теперь настало время вернуть долг. Докажем следующие утверждения:
	
	Предложение \ref{Fin->Periodic}: если конечная система запретов $S$ задает какое-то слово $w\in I$, то оно периодическое.\\
	\begin{Proof} 
		Пусть $k+1$ --- наибольшая длинна запрета в $S$. Рассмотрим $G_k^w$. Если это цикл --- то слово $w$ периодическое, если нет, то помимо наименьшего замкнутого пути $l$ соответствующего $w$ в этом графе есть еще хотя бы 1: поскольку $l$ самопересекается, можно рассмотреть замкнутый подпуть $l$ от точки самопересечения до нее же (меньший чем $l$ --- иначе это не точка самопересечения) --- и ходить по нему до бесконечности. Слово $w_1$ соответствующее этой циркуляции не равно $w$ (т.к. в нем встречаются не все слова длины $k+1$, которые встречаются в $w$). С другой стороны --- оно удовлетворяет $S$, т.к. содержит только те подслова длины $k+1$, которые встречаются в $w$, и если какое-то слово из $S$ содержится в $w_1$, то оно содержится в каком-то его подслове длины $k+1$ --- а значит и в $w$. Противоречие.
		
		\includegraphics[height=23mm]{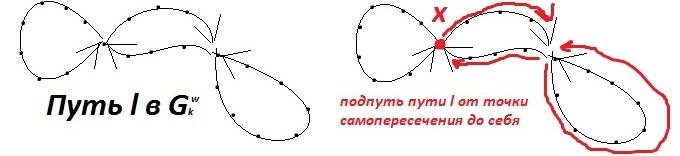}
		
		Заметим, что этим мы доказали еще один общий факт: любое бесконечное периодическое слово, соответствующее некоторому замкнутому пути в $G_k^w$ удовлетворяет $S\cap F_{k+1}$ (если $S$ --- приведенная система запретов задающая слово $w$).
	\end{Proof}

	Предложение \ref{Periodic->existsFin}: для любого $w\in I_n$ существует конечная система запретов, задающая его. \\
	\begin{Proof}	
		Расмотрим $n$-й граф Рози $G_n^w$ этого слова. Мы доказали, что это цикл длины $n$. Очевидно, слово $w$ удовлетворяет системе запретов $S=S_w\cap F_{n+1}\subset S_w$. Если какое-то другое слово $w_1$ тоже удовлетворяет ей, то $G_n^{w_1}\subseteq G_n^w$, т.к. он не может содержать никаких других ребер. Чтобы в нем был хоть один бесконечный в обе стороны путь (соответствующий, собственно, слову $w_1$) необходимо $G_n^{w_1}=G_n^w$. Тогда этоть путь единственен, а слово $w_1$ соответствующее ему --- совпадает с $w$. Противоречие.
	\end{Proof}

	\begin{ssSt} 		
		\label{reductedoao}
		Для данного $w\in I_n$ существует единственная приведенная система запретов, задающая его, и состоит она в точности из тех слов, которые обнаруживаются при рассмотрении графов Рози.
	\end{ssSt}
	\begin{Proof} 
		Тривиально следует из предложений \ref{Periodic->existsFin}, \ref{Reduction}  и \ref{obligatoryRestrict}
	\end{Proof}
	
	Заметим, что в этом разделе, как и в разделе 2.1 не используется двухсимвольность алфавита $A$
	
\subsubsection{Основной вывод}

	Из лемм \ref{forkDiff} и \ref{graphFin}, а также предложения \ref{reductedoao} очевидно
	\begin{Foll} 		
		\label{reductNum}
		$|S| = \sum\limits_{k=1}^\infty{|S\cap F_k|} = \sum\limits_{n=1}^\infty{(d_{k-1}^w - d_k^w + c_{k-1}^w)} = \sum\limits_{k=0}^\infty{c_k} + d_0$
	\end{Foll}

	Начинать рассмотрение целесообразно с $G_0$ (вершина --- пустое слово, и два ребра --- "`$a$"' и "`$b$"'). При этом можно не рассматривать тривиальные случаи когда запрещено $a$ или $b$, т.к. тогда оценка выполняется и точна (период слова состоящего только из буквы "`a"' или "`b"' равен 1). Тогда начальное количество развилок равно 1. $G_0$:
	
	\includegraphics[height=10mm]{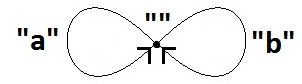}

\subsection{Размер графов Рози}

	Изучим как изменяется $|G_k|$ (количество ребер в $G_k$). 
	\begin{Lemm} 		
		\label{edgeDiff}
		$|G_{k+1}^w| = |G_k^w| + (d_k^w + c_k^w) - |S_w\cap F_{n+1}|$
	\end{Lemm}
	\begin{Proof} 
		Поскольку ребра $h(G_k)$ взаимно соответствуют путям длины два в $G_k$, их количества равны. Сопоставим каждому пути его первое ребро. Значит, количество таких путей равно $|G_k|$ плюс количество ребер, конечная вершина которых имеет исходящюю степень 2. Количество таких ребер равно $2c_k + f_k = t_k = d_k + c_k$. \parskip=3pt
	
		\noindent $|G_{k+1}^w| - |h(G_k)^w| = |S_w\cap F_{n+1}|$: на каждый запрет мы стираем по ребру. 
	\end{Proof}
	
	Значит, согласно леммам \ref{graphFin} и \ref{edgeDiff}, а также следствию \ref{reductNum} верно:
	\begin{Foll} 		
		$\forall w\in I_n$ верно равенство:
		\begin{multline*}
			n = |G_n^w| = |G_0| + \sum\limits_{k = 0}^{n-1}{(|G_{k+1}^w| - |G_k^w|)} = \\
			= |G_0| + \sum\limits_{k = 0}^{n-1}{((d_k^w + c_k^w) - |S_w\cap F_{n+1}|)} = \\
			= |G_0| + \sum\limits_{k = 0}^{n-1}{d_k^w} + \sum\limits_{k = 0}^{n-1}{c_k^w} - |S| = |G_0| + \sum\limits_{k = 1}^{n-1}{d_k^w}
		\end{multline*}
	\end{Foll}
	
	По сути это значит что на каждую $d$-развилку в каждом графе Рози приходится увеличение конечного количества ребер (и, соотвтетственно, периода) на 1. Остается только найти разумный способ сопоставить некоторые группы $d$-развилок запретам (или, поскольку их количества равны, перекресткам) так, чтобы получить необходимую оценку --- что мы в конечном счете и сделаем (на самом деле, эти группы --- те самые $md$-развилки, и каждая из них будет сопоставляться перекрестку, которым заканчивается, но как видно из утверждения теоремы рост конечного количества ребер от количества перекрестков может быть экспоненциальным, и значит важно в каком порядке 
	
	
\subsection{Другие способы нарисовать $G_{k+1}$ по $G_k$}

	Сейчас мы опишем и докажем корректность другого, необходимого для доказательства способа построения $G_{k+1}$ по $G_k$. Для лучшего понимания, упрощенного обоснования и ясной мотивировки мы приведем 8 разных способов построения, изменяющихся постепенно от имеющегося способа к необходимому.
	
	\vspace{-6pt} Для первого способа нам понадобится понятие схемы.
	\begin{ssDef}		
		Схемой ориентированного графа G будем называть взвешенный ориентированный граф $\widetilde{G}$ множеством вершин которого является множество вершин $G$ не являющихся дорогой, а множеством ребер --- пути в $G$ соединяющие такие вершины и не содержащие ни одной из них кроме как своим началом и концом (с весами равными длинам этих путей). Естественно, ориентация ребер в путях должна быть согласована и наследуется соответствующему ребру (из этого ясно, что такие пути не самопересекаются). \\[-20pt] 
	\end{ssDef}
	
	\includegraphics[height=38mm]{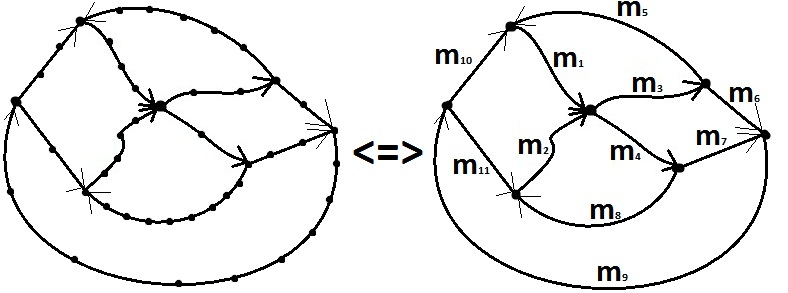}

	\vspace{-8pt} Другими словами, схема графа Рози --- это взешенный ориентированный граф, полученный из графа Рози заменой цепочек идущих подряд без развилок ребер на одно ребро с весом равным длине этой цепочки (и соответствующей ориентацией). Она наглядно отражает структуру графов.
	
	\vspace{-8pt} Формально: рассмотрим отношение эквивалентности порожденное следующими равенствами: два последовательных ребра эквивалентны, если их общая вершина --- дорога. Очевидно, в каждом классе эквивалентности (если только граф не является циклом) будет ровно по одному ребру не начинающемуся и не заканчивающемуся в "`дороге" ($b$ и $e$). Построим новый взвешенный ориентированный граф с множеством вершин равным множетсву развилок и перекрестков, а множеством ребер --- для каждого класса эквивалентности по ребру: из начальной вершины ребра $b$ в конечную вершину ребра $e$ с весом равным количеству ребер в классе эквивалентности. Также легко производится обратная операция (при этом вершины, соединенные ребром с весом 0 отождествляются). 
	
	\vspace{-8pt} Заметим, что на схемы корректно переносятся все определения развилок --- вес можно просто не учитывать. 
	
	\vspace{-8pt} Мотивация введения схем: это необходимо для введения фиктивных ребер с весом ноль в некоторых местах с целью удобства подсчета периода равного сумме приростов количеств ребер, соответствующих $md$-развилкам, соответствующих запретам.
	
	\vspace{-8pt} Используя понятие схемы предложим первый новый способ нарисовать $G_{k+1}^w$ по $G_k^w$: 	
	\begin{Meth} 
		\label{Basic} 
		Рассмотрим $\widetilde{G_k^w}$. Заменим в нем перекрестки на фиктивную конструкцию эквивалентную той, что получается из перекрестка при полушаге, только с весами 0 (см. рисунок) --- назовем эту операцию 0-заменой. Тогда перекрестков не останется --- только развилки. $\widetilde{G_k^w}$ изменим по правилу: если ребро закначивается однотипными развилками, его вес не изменится, если разными --- то $+1$, если его начало --- входящая развилка и $-1$, если исходящая (по сути --- это тот же полушаг). Теперь по схеме восстановим обратно граф (заменой ребер с весами на цепочки --- и получится точно такой же граф, как после полушага), а после этого применим запреты. 
	\end{Meth}
	\includegraphics[height=28mm]{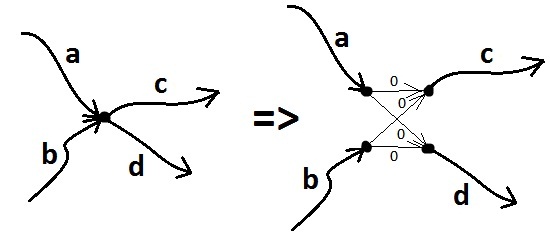}
	\begin{Proof} 
		Для доказательства корректности этого метода достаточно доказать что после перехода обратно от схемы к графу мы получим такой же граф, как после полушага. Заметим, что количество $t$-развилок в $G_k$ очевидно равно количеству развилок в схеме после 0-замены, и они естественным образом отождествляются: если вершина $t$-развилки была развилкой, то ей она и соответствует, а если перекрестком, то ей соответствует та развилка, в которой содержится ее "`направляющее"' ребро (оно принадлежит ровно одному из путей a,b,c,d --- см. рисунок). Каждой $t$-развилке соответствует $d$-развилка в $h(G_k)$. Остается заметить, что длины путей между $d$-развилками в $h(G_k)$ --- точно такие же как веса ребер между соответствующими им развилками в нашей схеме после преобразования --- возьмем пару развилок в $\widetilde{G_k}$ соединенных некоторым ребром. Если его вес 0, то после преобразования он будет равен 1 (т.к. это могли быть только развилки получившиеся заменой перекрестка) --- как и после полушага, а если нет, то см. рисунок.
	\end{Proof}
	\includegraphics[height=35mm]{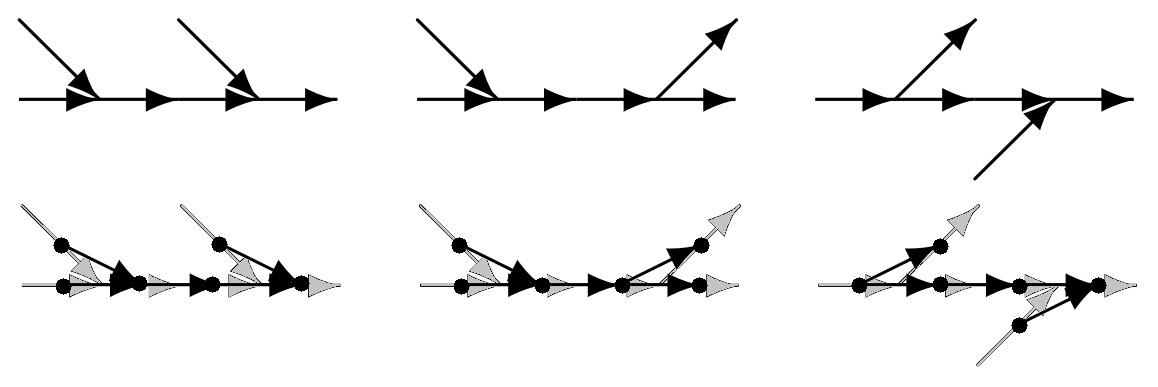}\label{picture_3_distance_transformation}
	
	\begin{Des}{
		$h^{*}(G)$ --- результат применения к графу G такой операции (пути между однотипными развилками не меняем, а между разнотипными --- сокращаем на 1 если начало --- входящая, а конец --- исходящая, и увеличиваем на 1, если наоборот).}
	\end{Des}
	
	Заметим, что метод \ref{Basic} принципиально отличается от метода \ref{Trivial} тем, что описывает преобразование графа не локально, а глобально: становится видно, что структура графа меняется только в местах перекрестков, а в остальных местах только изменяются расстояния. Причем противоположные развилки как бы движутся навстречу друг другу в направлении стрелок которые изображают. Заметим также, что перекрестки образуются по сути при встречах противоположных развилок.

%
	\begin{Meth} 
		Второй способ отличается от первого только тем, что запреты применяются до перехода от схемы обратно к графу Рози (удаляется соответствующее ребро).
	\end{Meth}
	\begin{Proof}
		Для доказательства корректности этого метода вспомним, что как уже было объяснено в доказательстве предложения \ref{restrictDiff} запреты происходят с ребрами для которых исходящая степень начальной вершины и входящая степень конечной равны 2. Но это значит что до перехода обратно от схемы к графу начальной и конечной вершинам соответствовали исходящая и входящая развилка, а нашему ребру --- ребро с весом 1 (и именно его мы удалим в нашей операции). Более того, это значит что до преобразования вес этого ребра был 0, и оно появилось в схеме после 0-замены.
	\end{Proof}

	\begin{Meth} 
		\label{Advanced}
		А третий --- тем, что теперь запреты вообще применяются до преобразования схемы. Результат будет эквивалентен, поскольку запреты ведь происходят только сразу после размножения, то есть в местах перекрестков --- а значит их можно сделать в самих наших фиктивных конструкциях. 
	\end{Meth}
	\begin{Proof}
		Поскольку уже показано что запреты будут применяться только к ребрам с весом 0 добавленным при 0-замене перекрестков, необходимо доказать только что если сначала применить запреты, а потом сделать преобразование схемы --- получится то же самое, как если сначала сделать преобразование, а потом применить запреты. Для этого достаточно доказать что графы будут иметь одинаковую структуру и соответствующие веса будут равны. Один случай, для наглядности, разобран на картинке. Поскольку структура схемы не меняется при нашем переходе (только веса) --- утверждение про идентичность структур очевидно (будут зачеркнуты одни и те же ребра). Заметим, что согласно нашей операции сумма весов на любом пути от одной развилки до другой при переходе меняется ровно так, как если бы на этом пути было одно ребро (очевидно из рассмотрения типов развилок на этом пути или картинки в доказательстве к первому методу --- путь между развилками там мог быть любой, не обязательно не содержащий других развилок) --- значит расстояния между оставшимися после наших действий развилками не зависят от последовательности этих действий.
	\end{Proof}
	Вообще, весь этот переход к схемам не обязателен, и для решения важно только иметь возможность разрешать перекрестки заменой данной врешины фиктивной структурой размера 0, но описывать преобразование графа "`там где развилки разнотипные одно ребро с весом один исчезло, а в другом месте ребро появилось"' все-таки менее удобно чем говорить про изменение веса.
	
	\begin{Meth} 
		Теперь все-таки попробуем обойтись без перехода к схемам. Просто изменим каждый путь между двумя развилками в графе соответствующим образом (полушаг), а потом применим запреты. При этом перекрестки можно обрабатывать по-разному: можно считать что у всех ребер бывших в $G_k$ веса 1, а с перекрестками делать 0-замену, а затем делать описанное только что преобразованее (при этом ребро веса 0 --- как бы путь длины 0 --- заменяется на 1 ребро с весом 1 --- путь длины 1). Или можно считать что в месте перекрестка находятся 4 $t$-развилки на расстоянии 0 (между противоположными из них), и когда рассматривается какой-то путь содержащий своим концом перекресток --- он считается путем до той развилки, направляющее ребро которой в нем содержится). 
	\end{Meth}
	\begin{Proof}
		Этот метод очевидно эквивалентен методу \ref{Basic}.
	\end{Proof}
	
	\begin{Meth} 
		Теперь сначала полностью разрешим перекрестки, а потом сделаем полушаг.
	\end{Meth}
	\begin{Proof}
		Подробнее: Рассмотрим граф Рози $G_k$. Каждому его ребру присвоим вес равный 1. Произведем 0-замену. Удалим ребра, соответствующие ребрам удаляемым в $h(G_k)$ (мы знаем, что удаляемые из $h(G_k)$ ребра соответствуют ребрам в схеме $\tilde{h(G_k)}$ с весом 1 и противоположными развилками на концах, а они соответствуют ребрам с весом 0 в текущем графе). Изменим длины путей в графе по описанным ранее правилам и отождествим вершины соединенные ребрами с весом 0. Веса забудем. Получим $G_{k+1}$ --- это следует из метода \ref{Advanced}.
	\end{Proof}

	Согласно замечанию о том, что длина любого пути между двумя развилками меняется только в зависимости от их типов, можно представить себе, что, согласно тому что расстояния между однотипными развилками не меняются, одни из них стоят на месте, а другие --- едут им на встречу со скоростью один.
	
	\begin{Meth} 
		С одной стороны --- этот способ принципиально отличается от всех предыдущих, и содержит существенно важную для решения идею, а с другой --- отличается от предыдущих только взглядом на вещи. Суть его вот в чем: сначала, как уже много раз делали, заменим все перекрестки на фиктивные конструкции размера 0 (очевидно, перекрестков у нас после этого не останется). А затем сделаем некоторую новую операцию, результатом которой будет то же самое что при полушаге, но процедура другая: выберем тип развилок (входящие или исходящие). Предположим без о.о. что входящие. Будем считать, что входящие $t$-развилки стоят на месте, а исходящие --- сдвигаются на 1 в ту сторону в которую они "`указывают"', поглощая одно ребро из пути перед собой, и выплевывая два --- по одному в каждый из хвостов за собой. 
	\end{Meth}
	
	\begin{picture}(50,60)
\put(0,30){\vector(1,1){10}}
\put(0,50){\vector(1,-1){10}}
\put(10,40){\vector(1,0){10}}
\put(20,40){\vector(1,0){10}}
\put(30,40){\vector(1,1){10}}
\put(30,40){\vector(1,-1){10}}
\put(40,50){\vector(1,0){10}}
\put(40,30){\vector(1,0){10}}
\put(0,0){\vector(1,1){10}}
\put(0,20){\vector(1,-1){10}}
\put(10,10){\vector(1,0){10}}
\put(20,10){\vector(1,1){10}}
\put(20,10){\vector(1,-1){10}}
\put(30,20){\vector(1,0){10}}
\put(30,0){\vector(1,0){10}}
\put(40,20){\vector(1,0){10}}
\put(40,0){\vector(1,0){10}}
\end{picture}

	\begin{Proof}
		Заметим, что все пути между развилками изменятся строго по нашему правилу. Значит этот метод очевидно эквивалентен методам \ref{Basic}, 2 и 4.
	\end{Proof}
	
	\begin{Meth} 
		То же самое, но теперь будем сначала сразу разрешать все перекрестки, а потом сдвигать исходящие развилки к входящим.
	\end{Meth}
	\begin{Proof}
		Эквивалентен методам 3 и 5.
	\end{Proof}
	
	\begin{Meth} 
	\label{parallel_drive}
	Строго говоря, это уже не способ получить следующий граф Рози по предыдущему, и графы Рози данного слова будут появляться только как промежуточная фаза каждой такой операции, но поскольку для исследуемого нами вопроса важно только количество произведенных за все операции запретов и количество ребер в конечном цикле --- это не важно. Здесь мы будем сначала делать полушаг в новом смысле, а потом разрешать все образовавшиеся перекрестки. Начинать в этом случае будем с графа $\widehat{G_0} = \mathbb G_2$ полученного из $G_0$ разрешением перекрестка. Варианты $\mathbb G_2$ (такой выбор индекса будет пояснен чуть позднее):
	\end{Meth}
	\includegraphics[height=20mm]{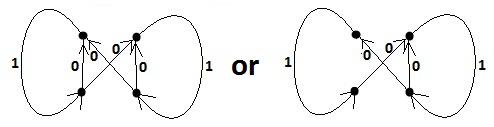}
	\begin{Proof}
		Однако нас интересует только общее количество запретов (или, эквивалентно, общее количество перекрестков) встретившихся в графах в течение всего процесса и (так как очевидно, как и графы Рози наши графы станут с некоторого момента одинаковыми циклами) --- конечное количество ребер. Они, очевидно, будут равны соответствующим количествам для последовательности графов Рози некоторого данного слова (поскольку последовательность, фактически, выглядит так: $\widehat{G_0} \rightarrow G_1\rightarrow\widehat{G_1}\rightarrow G_2\rightarrow\cdots$
	\end{Proof}

\subsection{Что мы получили}

	Имеется некоторое бесконечное в обе стороны периодическое слово $w \in I_n$ и задающая его приведенная система запретов $S$. Имеется последовательность графов Рози этого слова $G^w_0, G^w_1,\cdots, G^w_n$ (заканчивающаяся циклическим графом длина которого равна периоду слова n) и совокупность вспомогательных графов $\{h(G^w_i)\}$, с естественными биециями между множествами $E_{G^w_i} \leftrightarrow V_{h(G^w_i)}$ и $L^2_{G^w_i} \leftrightarrow E_{h(G^w_i)}$ ($L^2_{G^w_i}$ --- пути длины 2 в $G^w_i$, а также множествами $V_{h(G^w_i)} \leftrightarrow V_{G^w_{i+1}}$ и $E_{h(G^w_i)} \leftrightarrow E_{G^w_{i+1}} \sqcup (S \cap F_{i+1})$
	
	Также имеется последоваетельность графов $\widehat{G^w_0}, \widehat{G^w_1},\cdots, \widehat{G^w_n}$, таких что $h^{*}(\widehat{G^w_i}) = G^w_{i+1}$, которые строятся индуктивно по следующим правилам: $\widehat{G^w_i} \rightarrow h^{*}(\widehat{G^w_i}) \stackrel{0-exchange}{\longrightarrow}\widehat{G^w_i}'\stackrel{restrictions}{\longrightarrow}\widehat{G^w_{i+1}}$, где запреты происходят в тех же ребрах, где и в $h(G^w_i)$, структурно равном $\widehat{G^w_i}'$.
	
\subsection{Что дальше}

	В процессе применения 8-го метода абсолютно наглядно на каждую из $d_k^w$ $d$-развилок какого-то типа за 1 шаг (от $\widehat{G^w_{k-1}}$ к $\widehat{G^w_k}$) приходится прирост общего количества ребер на 1, на каждый перекресток (из $c_{k+1}^w$) --- прирост количества $d$-развилок на 1, на каждый запрет --- уменьшение количества $d$-развилок на 1, и все эти приросты-уменьшения происходят параллельно. Ключевая идея доказательства оценки --- в том чтобы делать это последовательно, на каждом шагу с какой-то одной развилкой, с которой это сделать возможно: если в данном графе перед какой-то исходящей (в направлении ее движения) развилкой до ближайшей входящей нет других исходящих --- единомоментно проведем ее до этой входящей развилки через этот путь и разрешим образовавшийся перекресток (применив запреты так же, как мы собирались это сделать при планомарном изменении графа --- параллеельном движении развилок). Заметим, что поскольку количество запретов равно количеству перекрестков --- мы можем считать только последние. Так, при каждой такой операци количество перекрестков увеличивается на 1, а $|\mathbb G_k|$ увеличивается на длину пути, по которому мы провели развилку (вот мы и сопоставили каждому запрету прирост количества ребер --- ту самую $md$-развилку). Соответственно целесообразно индексировать наши графы количеством учтенных перекрестков (поэтому и индекс начального графа --- 2: для превращения $\mathbb G_2$ в цикл необходимо два запрета). Как уже можно было заметить, термин "`шаг"' используется автором для любой минимальной циклически повторяющейся операции над графом. Теперь мы будем использовать его для обозначения последней из описанных --- передвижения развилки и разрешения появившегося перекрестка.
	
	Теперь фиксируем слово ($w \in I_n$) и его приведенную систему запретов $S$ (впрочем, и так единственную).
	
	Рассмотрим $\mathbb{G}_2 = \widehat{G_0}$. Как мы уже говорили, есть два, с точностью до перемены $a$ и $b$ местами варианта. Рассмотрим в последовательности $\widehat{G_i}$ $l$-развилки (развилки естественным образом отождествляются с наследниками в методе \ref{parallel_drive} --- если не произошло 0-замены). Вне завсисмости от того как выглядит $\widehat{G_0}$ для каждой развилки в нем каждому (из двух) ее хвостов (входящему в нее ребру для входящей или исходящему для исходящей) сопсотавим метку "`left"' или "`right"' (одну одному, а другую --- другому). Опишем как наследуются эти метки: при $h^{*}$ --- естественным образом (последнему ребру в пути, в котором было ребро "`left"' --- метку "`left"', аналогично для "`right"'). При 0-замене --- тоже довольно естественно:
	
	\includegraphics[height=35mm]{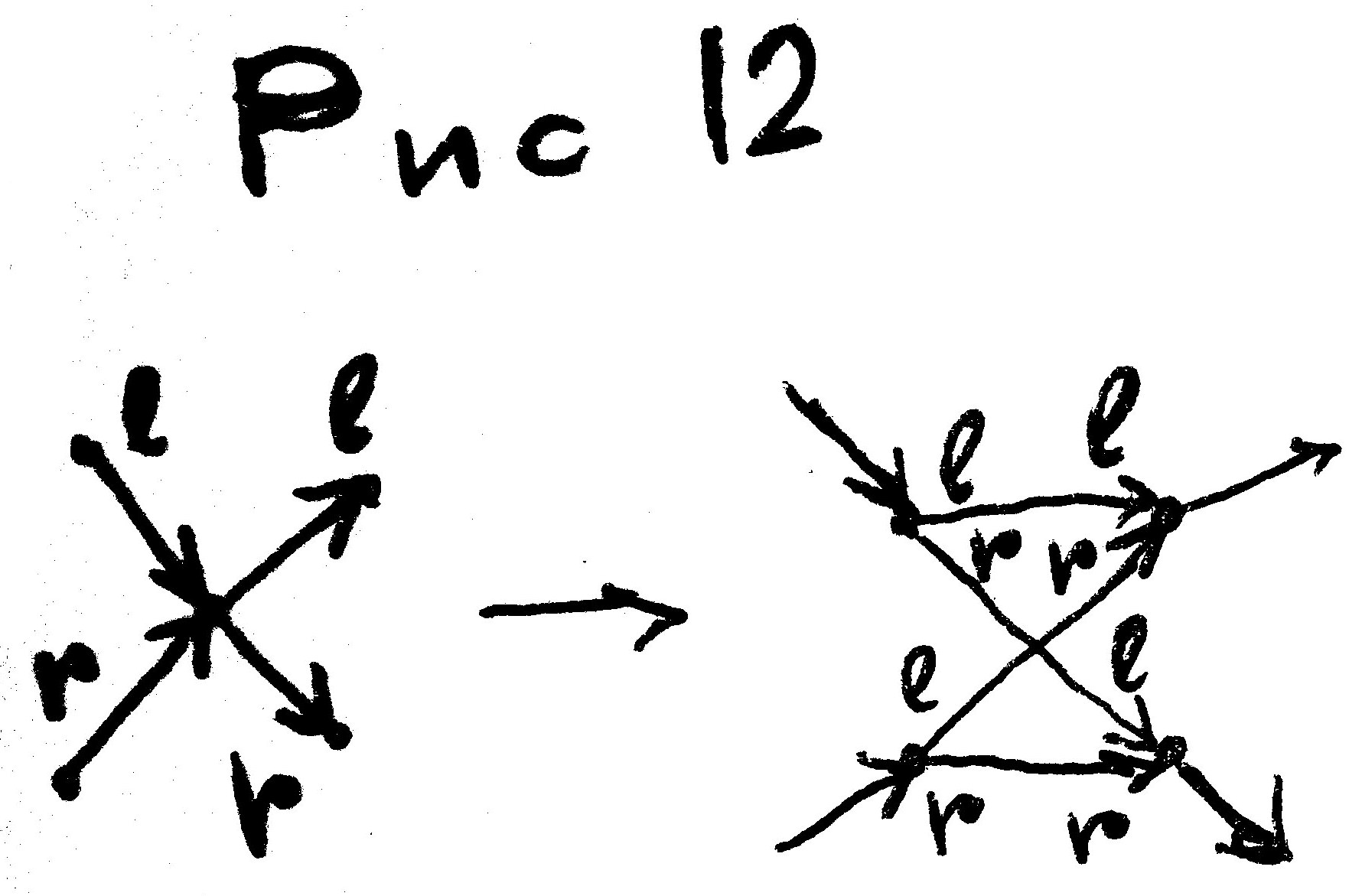} 
	
	(После применения запретов какие-то из нарисованных развилок исчезнут, а какие-то останутся. Вместе с метками). В общем-то получается, что на время существования $l$-развилки в графе для нее вполне естественно определены левый и правый хвост. Понятно, что с помощью этих меток однозначно определяется какое ребро запрещается после 0-замены (даже --- меньшей информацией --- достаточно пометить развилки только одного типа --- допустим, только входящие). Также заметим, что для каждой $l$-развилки однозначно определяется ее история (после взаимодействия каких пар развилок противоположного типа она появилась и левым или правым хвостом являлись соответствующие развилки) и то с какой развилкой она в итоге размножилась (из множества всех $l$-развилок противоположного типа), да и вообще --- определена полная история взаимодействий (какие пары развилок размножились, какие ребра были запрещены и кто были их дети).
	
	Возьмем $\mathbb{G}_2$. Расставим в нем как-нибудь метки (и установим соответствие с $l$-развилками в $\widehat{G_0}$). Будем производить в произвольном порядке следующие (уже ранее описывавшиеся) операции: если в текущем графе на пути между двумя данными развилками (входящей и исходящей --- в указанном порядке) нет других развилок (причем путь проходит через их направляющие ребра), то сдвинем конкретно эту исходящую (для определенности) развилку по описанным правилам (одно ребро поглощает, два выдает) до входящей, произведем 0-замену образовавшегося перекрестка, расставим метки "`left"' и "`right"', удалим соответствующие ребра (такие же, как те, которые были удалены в изначальной последовательности графов $\widehat{G_i}$), оставшимся развилкам поставим в соответствие $l$-развилки этой последовательности --- те, которые получились при размножении соответствовавших текущим развилкам $l$-развилок. 
	
	\begin{sSt} 
	В таком случае эта операция всегда определена корректно (те $l$-развилки которые соответствуют двум развилкам с которыми мы производим данную операцию в исходной последовательности размножались именно друг с другом), в конце эта последовательность станет циклом (эта часть очевидно следует из предыдущей) и длина цикла будет равна длине финального цикла исходной последовательности графов --- вне зависимости, в частности, от выбора на каждом шаге конкретной пары развилок между которыми ничего нет. 
	\end{sSt}
	\begin{Proof}
		Вообще, это очевидно --- если представить, что в исходной последовательности развилки двигаются паралельно (не обгоняя друг друга), встречаются, размножаются, происходят какие-то запреты --- то, если все происходит идентично --- очевидно не важно в какой последовательности это происходит если они друг друга не обгоняют (т.к. при встрече развилок пути переклеиваются идентично).\\[-7pt]
		
		Рассмотрим схемы графов из процесса применения метода \ref{parallel_drive}. Как мы уже говорили, при переходе от графа к графу меняются только веса ребер и небольшие фрагменты графа в местах встречи развилок (когда ребро становится равным 0). Заметим, что согласно с этим ребра могут появляться, исчезать и объединяться. Всё это --- только при встрече развилок. Рассмотрим множество всех ребер когда-либо появившихся в этих графах $\mathbb{E}$ (можно считать что при встрече развилок появляются 4 ребра, а потом мы сразу некоторые из них убиваем, а можно --- что появляются только те из них которые остаются, и некоторые сразу объединяются с теми между которыми убита развилка. Для удобства и ясности примем первое) . Мы помним, что у нас уже есть множество размеченных развилок, и соответственно легко установить в последовательности $\mathbb{G}_k$ какие ребра каким соответствуют при появлении (опять же рассматривая схемы). При объединении двух ребер будем писать на нем слово в алфавите $\mathbb{E}$ получающееся приписыванием текущего слова на первом ребре к текущему слову на втором (если объединяются три подряд --- тоже понятно как). При появлении на ребре будем писать однобуквенное слово из того ребра что ему соответствует. Заметим, что в каждом графе на концах ребра находятся какие-то развилки. Если развилка "`смотрит"' в сторону ребра, концом которого она является --- она, очевидно, будет его концом до исчезновения этого ребра (возможно после его объединения с какими-нибудь другими ребрами на другом конце) из-за встречи этой развилки с противоположной. Если эта развилка смотрит "`от"' --- она встретится с кем-нибудь там с той стороны, и у этого ребра либо поменяется конец на развилку смотрящую "`к"', либо припишется с этой стороны какое-нибудь ребро. Утверждение: в каждом графе $\mathbb{G}_k$ на ребрах будут написаны только слова, являющиеся подсловами слов написанных на ребрах процесса метода \ref{parallel_drive} в момент их смерти (то есть после максимального объединения слов), а также у каждого ребра на конце будет та развилка, которая должна быть --- либо смотрящая в его сторону и одновременно та, с которой он должен умереть (можно считать что развилки тоже буквы --- тогда просто это то самое условие на подслова), либо смотрящая в противоположную сторону и та, которая при размножении заменится соотв. на нужную смотрящую в его сторону, либо приведет к приклеиванию однобуквенного слова --- того, которое идет перед данным словом (в том самом длинном слове в котором оно встречается как подслово --- очевидно такое одно, т.к. буквы во всех словах разные и не повторяются).\\[-7pt]
		
	Заметим, что по сути наше утверждение (о том что операции можно делать в любом порядке) --- есть утверждение об ассоциативности произведения операций. Соответственно его доказательство также будет сродни классическим доказательствам ассоциативности.  \\[-7pt]
	
	Всё это доказывается по индукции. Для $\mathbb{G}_2$ это верно, т.к. там всё так же как в $G_1$. Пусть для $\mathbb{G}_k$ данное утверждение верно. При переходе к $\mathbb{G}_{k+1}$ мы производим всего одну операцию --- встречаются две развилки находящиеся на концах ребра, с которым они должны умереть. Очевидно значит, они размножались именно друг с другом в основном процессе (т.к. вместе с одним ребром умирают только две конкретные развилки). Нас интересуют четыре ребра концами которых являлись эти развилки с других сторон. Условия на граф обеспечивают то что при 0-замене к уже имевшимся ребрам приклеились именно те развилки или ребра, которые нужно. Остается только проверить, что для четырех новопоявившихся ребер верно следующее: либо они исчезнут, либо у них на концах те развилки которые нужно. Это просто: поскольку появление ребер и развилок естественным образом привязано друг к другу (мы знаем что в основной последовательности графов эти ребра и развилки появились вместе --- мы их так сопоставляем) --- то на концах ребер развилки смотрящие "`от"', после которых прибавится то что нужно --- либо развилка с которой мы встретим сметрь, либо однобуквенное слово, которое и должно быть (потому что так было в основном процессе). Корректность встречи как подслов индуцируется тем, что слова в основном процессе приписывались с "`таким"' упорядоченьем букв --- хоть может быть и в разной очередности во времени.\\[-7pt]
		
	Конечное количество ребер всегда будет одно и то же, т.к. оно просто равно сумме начального количества ребер плюс для каждого ребра на котором написано финальное слово (то, с которым это ребро умрет, а не к нему что-нибудь приклеится) --- количество раз сколько его левая граница сдвигалась вправо, и количество раз сколько его правая граница сдвигалась влево (в нашем предположении одни из них всегда стоят).\\[-7pt]
	
	Будем считать, что все наши действия разбиты на малюсенькие части --- сдвиги одной развилки на 1. При этом происходит -1 к длине одного ребра и +1 к длинам двух. Ребра делятся на те, которые когда-нибудь умрут, и те которые не умрут (по сути --- будут частью слова написанного на конечном цикле графа). Заметим, что количество раз сколько концевые развилки для данного слова будут придвинуты друг к другу равно количеству раз сколько до этого развилки  в которых заканчиваются подслова этого слова будут отодвинуты (это почти очевидно). По сути --- конечная длина цикла --- количество раз, которые будут отодвинуты развилки на концах слов, которые не умрут. \\[-7pt] 
	
	Если развилки одного типа стоят, а другого --- двигаются, то нам не важно (для количества) в каком порядке они двигаются (если они не обгоняют друг друга) --- т.к. очевидно вне зависимости от порядка каждая развилка пройдет до своей пары одно и то же расстояние --- это можно доказать по индукции. Заметим, что любая развилка за время своего существования не меняет ближайшие к ней буквы написанные на ребрах заканчивающихся в ней. Ну, а поскольку всё корректно --- значит эти ребра являются подсловами одних и тех же конечных слов (т.к. данные буквы встречаются каждая только в одном конечном слове). Значит каждая развилка с момента появления до момента смерти дает прирост по +1 к словам содержащим как подслова оба ее хвоста. А время ее жизни равно сумме времен жизни всех тех развилок, которые давали прирост к ее слову (может быть --- удвоенным --- если оба хвоста содержались как подслово). По индукции по множеству развилок утверждение очевидно.		
	\end{Proof}

\section{Доказательство оценки. Техническая часть.}
\subsection{Утверждение оценки}

	Сформулируем утверждение оценки:
	
	 $|\mathbb G_k| \leq \dfrac{\varphi_{k+d_k+I_{(d_k>1)}-1}}{2^{d_k+I_{(d_k>1)}-1}}$, где $|\mathbb G_k|$ --- сумма весов всех ребер $\mathbb G_k$ (некорректно говорить об их количестве, при наличии ребер весом 0), $d_k$ --- количество входящих развилок развилок в $\mathbb G_k$, а $I_{(d_k>1)}$ равно 1, если $d_k>1$, и 0 иначе. Если $d_k = 1$ то дополнительно $2x + y + z \leq \varphi_{k+1}$.
	
	\includegraphics[height=20mm]{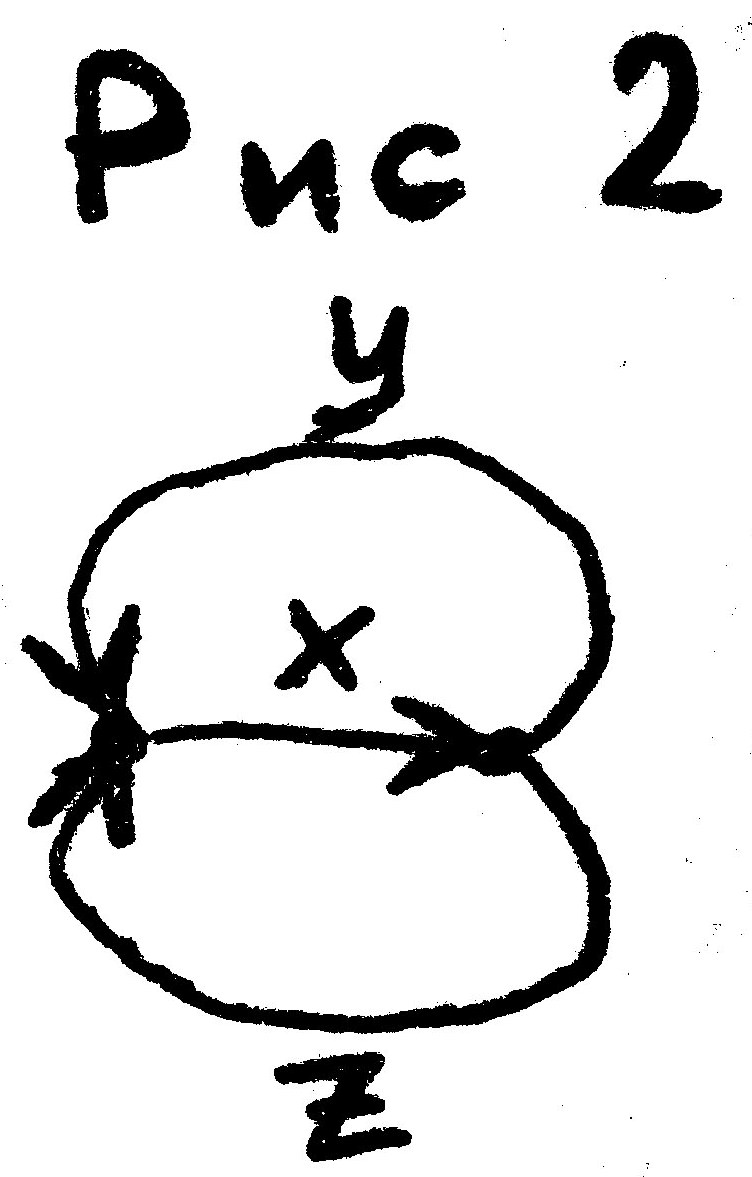}
	
	Заметим, что для графа $\mathbb{G}_1$ это утверждение верно:
	
	\includegraphics[height=20mm]{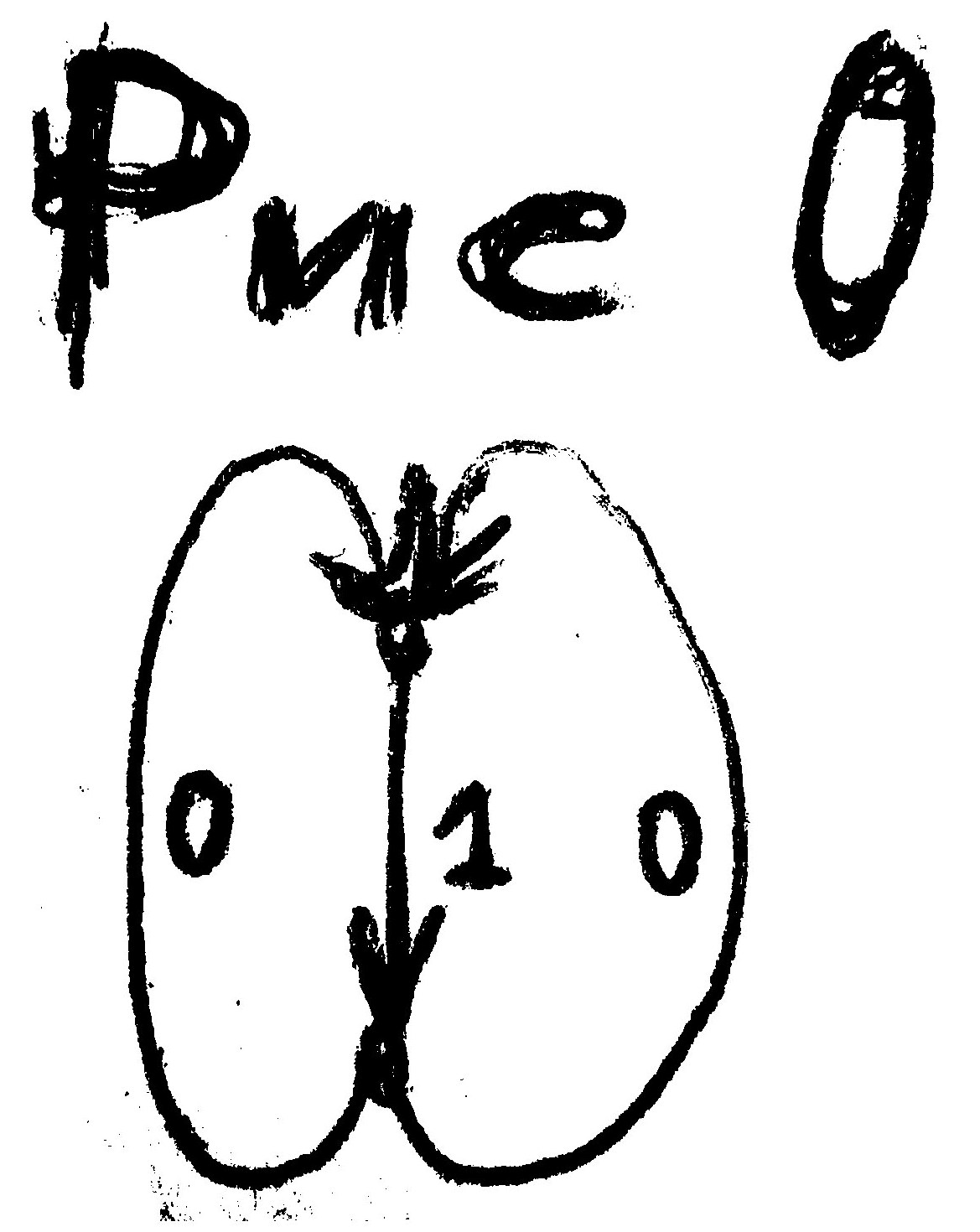}

	Мы будем начинать рассмотрение с этого графа и рассматривать его как базу индукции. Запреты при первом шаге будем делать так, чтобы получилось $\mathbb G_2$.
	
	Шаг индукции сформулируем следующим образом: если для некоторого $\mathbb G_k$ оценка верна, то можно сделать еще несколько шагов так, что для нового графа оценка также будет верна. Также будем утверждать, что для любого $\mathbb G_k$, такого что $d_k = 1$ --- оценка верна. 

\subsection{План доказательства}

	Доказательство естественным образом делится на три части. 
	Заметим, что если $d_k = 1$ --- у нас есть только один вариант какую пару развилок схлопнуть --- единственную имеющуюся в графе (под словом "`схлопнуть"' мы будем понимать произведение описанной ранее операции над данной парой развилок соединенных путем без других развилок). 
	Первая, самая простая часть доказательства касается случая, когда после этого единственного возможного схлопывания останется снова одна развилка ($d_{k+1} = 1$) --- в этой ситуации возможны только два случая, и в обоих оценка остается верной.
	Если $d_{k+1} > 1$, то оценка резко усиливается --- мы требуем запаса на одно гипотетическое удвоение. Однако, в силу того что начальная структура графа строго фиксирована, все возможные случаи перебираются, и либо оценка оказывается верна, либо мы возвращаемся к одноразвилочному состоянию и тогда она тоже верна.
	Третья часть доказательства соответственно касается случая, когда $d_k > 1$. В этом случае мы не знаем как выглядит граф, однако должны отследить только отклонение оценки в зависимости от количества произведенных операций и изменения количества развилок --- это позволяет обобщить большой класс ситуаций, когда мы можем сразу сделать несколько операций с сохранением верности оценки, а остальные опять разобрать перебором.

\subsection{Необходимые численные данные и мотивация}

	Мы будем изучать величину $\dfrac{|\mathbb{G}_j|}{|\mathbb{G}_i|}$ и оценивать ее некоторой константой $c$. Из верности утверждения оценки для $\mathbb{G}_i$ и неравенства $c \leq \dfrac{\varphi_{k + j - i + d_j - d_i}}{\varphi_k} * 2^{d_i - d_j}$ для $c$ (где $k$ --- соответствующий индекс для $\mathbb{G}_i$) будет следовать верность утверждения для $\mathbb{G}_j$. Если поделить и домножить правую часть этого неравенства на $\alpha^{j - i + d_j - d_i}$ где $\alpha$ --- золотое сечение ($\frac{1+\sqrt{5}}{2}$), то получится произведение двух выражений, значения которых легко анализировать: $\dfrac{\varphi_{k + j - i + d_j - d_i}}{\varphi_k * \alpha^{j - i + d_j - d_i}} * \left(\alpha^{j-i}\left(\dfrac{\alpha}{2}\right)^{d_j - d_i}\right)$ . Очевидно, что правая часть неравенства всегда будет больше чем второе выражение умноженное на минимальное значение первого. И почти всегда --- если умножить на 2-е по порядку значение первого (именно это значение нам будет удобнее всего, а оставшийся частный случай можно рассмотреть отдельно).
	
	Соответственно, приведем значения выражения $\dfrac{\varphi_u}{\varphi_v * \alpha^{u-v}}$ (Приложение \ref{table_evasion}, по вертикали откладывается $u$). Согласно нашей индексации $u,v \geq 1$. Нас интересует только нижний угол, где $u \geq v$, т.к. $j - i + d_j - d_i \geq 0$. Легко понять (из таблицы, а также из аналитической формулы для чисел Фибоначчи: $\varphi_c = \frac{ \alpha^{c+1}-(-\alpha)^{-c-1}}{\sqrt{5}}$), что 0.9270509831 (вычислено с б\'ольшей точностью) является наименьшим значением в нижнем углу таблицы, а 0.9442719100 --- вторым по порядку.
	
	Заметим также, что числа фибоначчи представляют из себя не чисто геометрическую прогрессию с показателем $\alpha$, а таковую с некоторыми небольшими (но существенными для асимптотики в начале) отклонениями, и таблица \ref{table_evasion} характеризует отклонение увеличение значения $\varphi$ от такового для геометрической прогрессии. 
	
	Значения выражения $\alpha^c*\left(\frac{\alpha}{2}\right)^k$ при $c>0$ и $-c\leq k\leq c + 1$, отражающего максимальное допустимое увеличение размера графа при $c$ шагах и изменении при них количества развилок на $k$ приведены в таблице \ref{table_pure} (в каждом элементе таблицы указано значение данного выражения, значение $c+1$ --- обычно, но не всегда в нашем доказательстве мы будем стремиться чтобы именно такое количество раз ограничивало увеличение размера графа за указанное количество шагов --- и значение $k$.)

	Приведем также таблицу значений этого выражения, домноженного на 0.9442719100 (небольшое количество случаев, когда $\dfrac{\varphi_u}{\varphi_v * \alpha^{u-v}}$ все-таки равно 0.9270509831 будет прокомментировано позднее): приложение \ref{table_second}. Первый столбец, где $k = c + 1 > c$ соответствует как раз начальному увеличению числа развилок (когда была 1, а стало больше) и понадобится в соответствующем доказательстве.

	Также заметим, что по данным таблицам (а также из очевидных соображений о производных) легко сделать выводы о росте соответствующих значений за их пределами.

\subsection{Одноразвилковая ситуация}

	Известно, что $d_k = d_{k+1} = 1, 2x_k + y_k + z_k \leq \varphi_{k+1}, x_k + y_k + z_k \leq \varphi_k$. Без ограничения общности возможны два случая: $x_{k+1} = x_k + y_k$, $y_{k+1} = x_k + z_k$, $z_{k+1} = 0$ или $x_{k+1} = x_k + z_k$, $y_{k+1} = x_k + y_k$, $z_{k+1} = 0$. В любом из них $|\mathbb G_{k+1}| = 2x_k + y_k + z_k \leq \varphi_{k+1}$, $2x_{k+1} + y_{k+1} + z_{k+1} \leq 3x_k + y_k + z_k + max(y_k,z_k) \leq (2x_k + y_k + z_k) + (x_k + y_k + z_k) \leq \varphi_{k+1} + \varphi_k = \varphi_{k+2}$
	
\subsection{Многоразвилковая ситуация}

	Заметим, что после любой операции количество развилок в графе может либо уменьшиться на 1, либо остаться тем же, либо увеличиться на 1. При этом размер графа увеличится на длину схлопнутого пути.
	
\subsubsection{Уменьшение количества развилок}

	Если количество развилок после некоторой операции уменьшилось на один, то оценка останется верна: очевидно, длина несамопересекающегося пути в графе не может превышать общего количества ребер в нем. Как мы уже отмечали, за одну операцию размер графа увеличивается на длину схлопнутого пути, то есть $|\mathbb{G}_{k+1}|  =  |\mathbb{G}_k| + s  \leq  2 * |\mathbb{G}_k|  =  2 * \dfrac{\varphi_{k+d_k+I_{(d_k>1)}-1}}{2^{d_k+I_{(d_k>1)}-1}}  =  \dfrac{\varphi_{(k+1)+d_{k+1}+1}}{2^{d_{k+1}}}$ (здесь очевидно $I_{(d_k>1)} = 1$). 
	
	Заметим, что если в графе останется одна развилка, то у нас будет запас на одно удвоение для асимптотики --- это не случайно.	Это используется для случая, когда $d_k$ было равно $2$, и, соответственно, $d_{k+1}$ стало равно $1$. Мы знаем, что $x_{k+1} + y_{k+1} + z_{k+1} \leq \dfrac{\varphi_{k+2}}{2} \leq \varphi_{k+1}$. Из этого следует, что $2x_{k+1} + y_{k+1} + z_{k+1} \leq 2(x_{k+1} + y_{k+1} + z_{k+1}) \leq \varphi_{k+2}$.
	
	Однако это еще не полное доказательство того, что для одноразвилкового графа оценка всегда верна. Дело в том, что у нас могло не быть промежуточной остановки на графе с больше чем одной развилкой с верной оценкой: если была одна развилка, потом стало больше (и усиленная оценка не выполнена), а потом обратно одна. Это случай будет рассмотрен в последней части доказательства.
	
	В силу простоты рассуждения в этом случае условимся, что если для некоторого графа оценка верна и в нем есть пара развилок, после схлопывания которой количество развилок уменьшается на 1 --- мы будем производить операцию с ней.
	
	Пусть теперь в графе нет таких пар развилок. 
	
\subsubsection{Когда пар развилок, которые можно схлопнуть много}
	
	Рассмотрим некоторый текущий граф для которого оценка верна. Схлопнем в нем все пары развилок которые есть (не схлопывая вновь образовавшиеся). В этих операциях размер графа увеличится не более чем в два раза (т.к. пути которые мы схлопнули и тем самым удвоили одновременно присутствовали в нашем начальном графе и не пересекались). Посмотрим на таблицу \ref{table_second}. То, что оценка перестала быть верна означает, что мы сделали 1 или 2 операции (причем если 2 --- то количество развилок увеличилось ровно на 2). Во всех остальных случаях оценка осталась верна. Это значит, что остаются только случаи когда в графе в паре находятся только 1 или 2 пары развилок. Их и будем рассматривать.
	
	Заметим, что вообще не быть схлапываемых пар развилок в графе не может: пойдем по ребрам от какой-нибудь входящей развилки в направлении ориентации. Пока на этом пути встречаются только входящие развилки путь однозначен. Идти всегда есть куда --- у нас нет вершин исходящей степени 0. Раз есть входящие, значит есть и исходящие. В силу конечности графа мы когда-нибудь в них придем --- иначе мы получим цикл, в который ребра только входят. Значит эти входящие развилки никогда бы не встретились с исходящими, и никогда не могли бы быть уничтожены --- противоречие, т.к. в конце мы получаем циклический граф. 
	
\subsubsection{Увеличение количества развилок при схлопывании}
	
	\textbf{Сначала рассмотрим случай когда пар развилок 2.} Схлопнем одну из развилок. 
	
	\includegraphics[height=30mm]{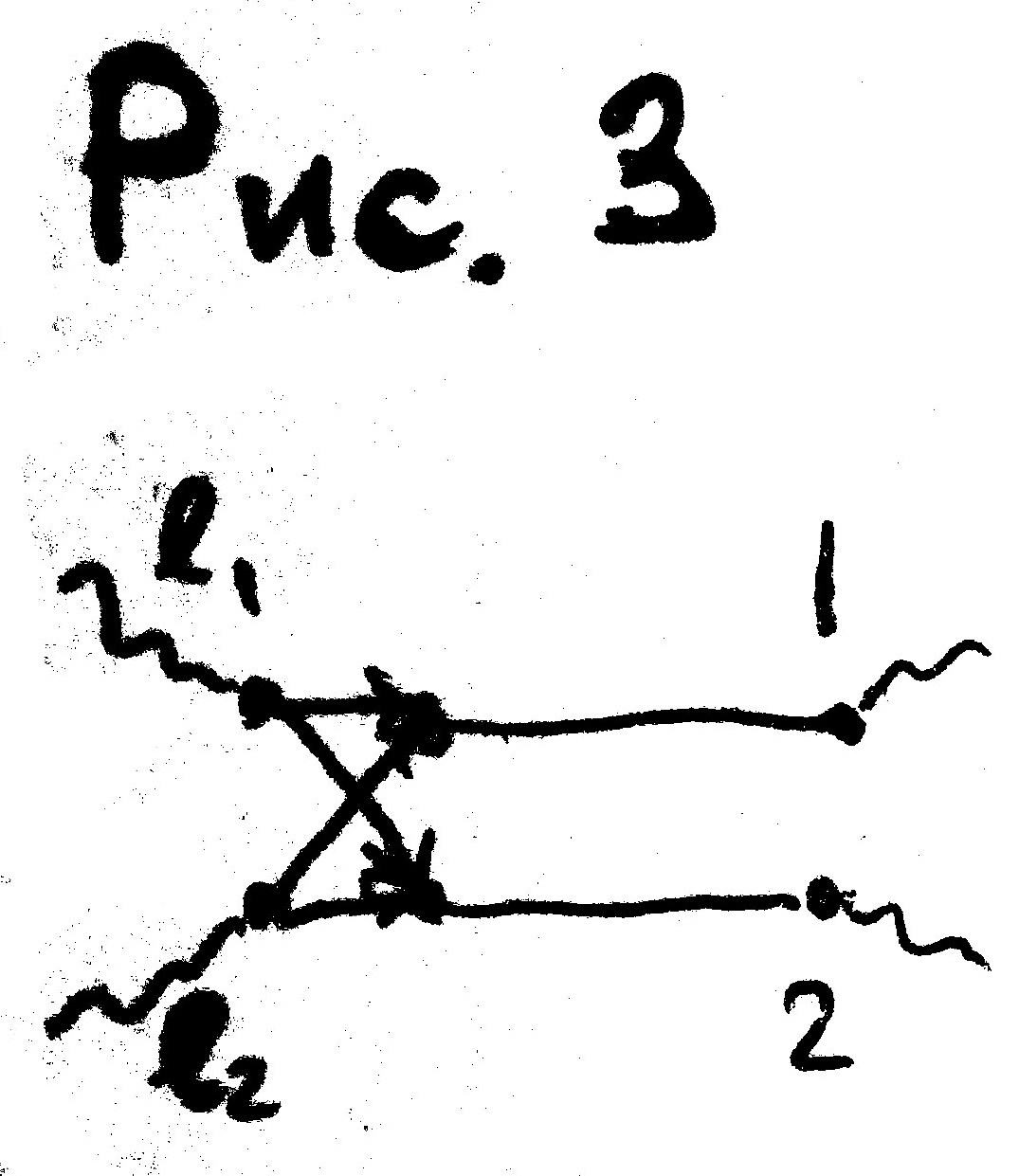}
	
	По одному из путей $l_1$ и $l_2$ точно есть какие-нибудь другие развилки (другой путь может непосредственно вести в вершину 1 или 2, но оба --- не могут, т.к. в графе другие развилки все-таки есть). Пойдем по этому пути пока не встретим входщую развилку. Если перед этим мы встречали исходящие развилки, то, очевидно, между этой входящей и последней исходящей ничего нет, и это та самая вторая пара, что есть в графе и что тоже схлопывается без запретов. Если нет --- то сейчас ничего нет между этой развилкой и развилкой в начале нашего пути. Во втором случае просто схлопнем эту пару и пару которая еще есть в графе из тех двух что были в начале. Получится, что граф не более чем удвоился. За три операции. Значит, исходя из нашей таблицы \ref{table_second} --- оценка осталась верна. В первом же случае --- сдалем так: в этой (второй) паре не исходящюю быстро проведем ко входящей, а наоборот. То, что так делать можно --- почти очевидно. Во-первых, потому что в самом начале мы могли решить остановить не входящие, а исходящие. Во-вторых, потому что структурно граф меняется точно так же, как если двигать исходящую ко входящей, а значит остается только понять корректность в смысле количественных характеристик --- количества ребер в конечном графе. Это мы докажем чуть позже. Но мы знаем, что до схлопнутой последней исходящей развилки на нашем пути по которому мы шли (в направлении противоположном ориентации ребер) были еще исходящие развилки (может быть, только одна --- в самом начале) Значит, поскольку мы не запрещали ребер при обоих схлопываниях, сейчас есть еще одна пара развилок. Схлопнем их. Три операции, не более чем удвоение (засчет того что мы удваивали только ребра присутствовавшие в графе до начала операций) --- оценка верна.
	
	\includegraphics[height=40mm]{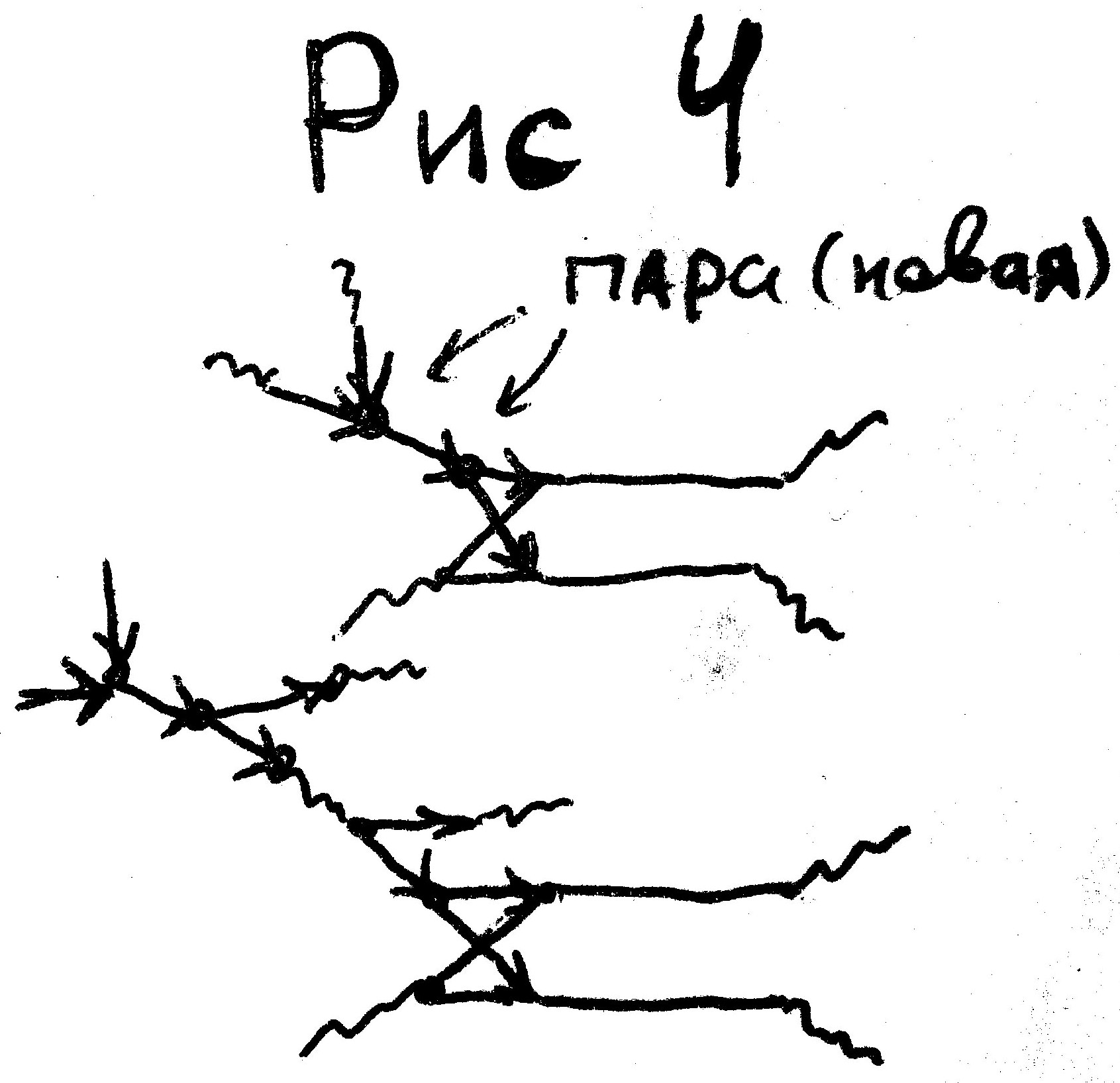}
	
	Докажем просто, что если в некоторый момент какую-то входящую развилку подвинуть на 1 вправо (по направлению ориентации ее направляющего ребра), а потом все делать как обычно (понятно, что так можно --- т.к. структуру графа это не поменяло) --- то конечное количество ребер будет тем же самым. Действительно, когда мы подвинем эту развилку, размер графа увеличится на 1. Но зато той развилке, которая схлопнется с этой будет идти на 1 меньше, и засчет этого мы потеряем 1. Когда эти развилки схлопнутся, останется сколько-то входящих и столько же исходящих развилок. Заметим, что теперь всем исходящим идти на 1 больше, чем в обычном процессе, а тем исходящим, которые схлопнутся с образовавшимися входящими --- на 1 меньше --- баланс опять 0. И так далее для каждого схлапывания, а всё остальное очевидно будет так же как обычно --- значит в конце ребер столько же.
	
	\textbf{Остается случай одной пары развилок.} В этом случае граф выглядит так: 
	
	\includegraphics[height=30mm]{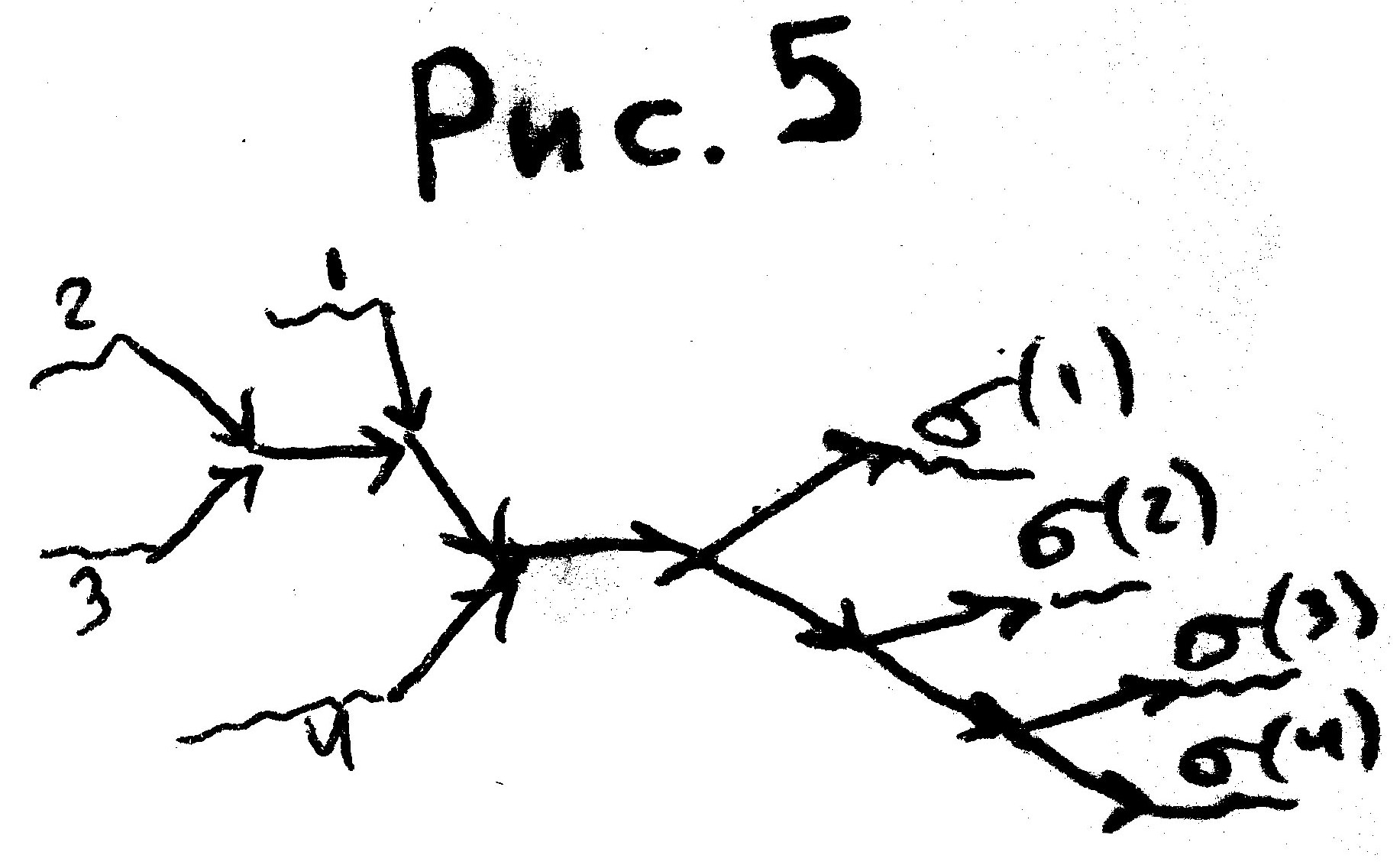}
	
	Пояснение --- слева за входящей развилкой некоторое дерево входящих развилок. Справа за исходящей --- дерево исходящих. Их листья соединены между собой в некотором порядке (перестановка $\sigma$).

	Соответственно мы можем только схлопнуть эту пару развилок. В этом пункте мы рассмотрим ситуацию, когда при этом не происходит запретов и количество развилок увеличивается. Тогда, поскольку в графе сейчас более одной развилки, хотя бы по одному из путей (аналогично предыдущим рисункам $l_1$ и $l_2$) будет входящая развилка. Схлопнем соответвующую пару. Две операции. Если во второй раз запреты были --- оценка верна (см. таблицу \ref{table_second}). Если нет --- мы, если в графе было более двух развилок, можем опять схлопнуть какую-нибудь вновь образовавшуюся по одному из путей пару развилок. Три операции. Не более чем удвоение. Оценка точно верна. Значит остается случай когда в графе ровно две развилки. 
	
	Без ограничения общности в этих случаях граф может иметь один из двух следующих видов (доказывается тривиальным перебором возможных перестановок с точностью до изоморфизма графов):
	
	\includegraphics[height=40mm]{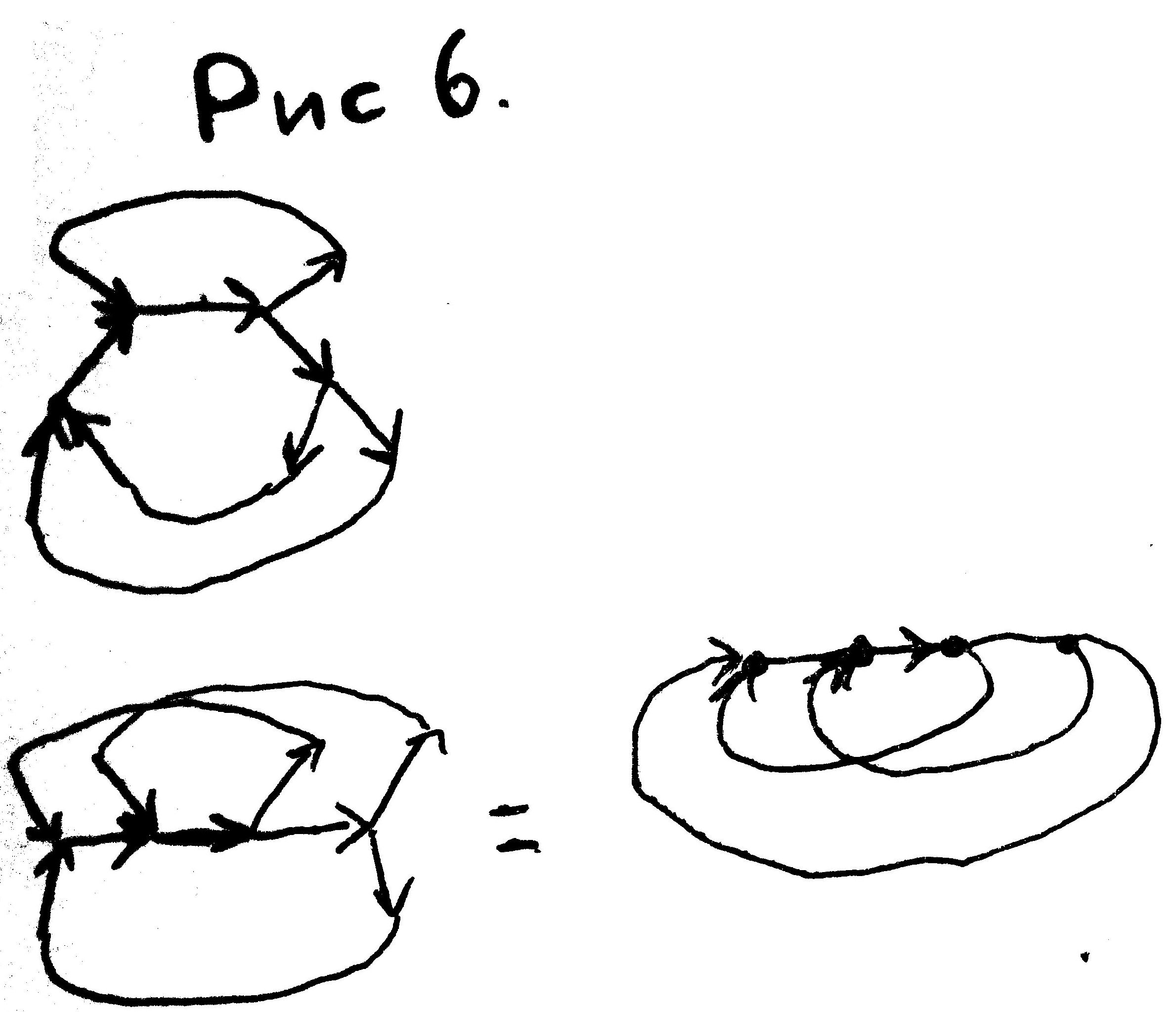}
	
	В обоих случаях начнем с того (с чего же еще?) что схлопнем единственную имеющуюся в графе пару развилок между которыми ничего нет. 
	
	Далее в первом случае, поскольку запретов по нашему предположению не произошло --- мы можем схлопнуть развилки 1, 2 и 3.
	
	\includegraphics[height=30mm]{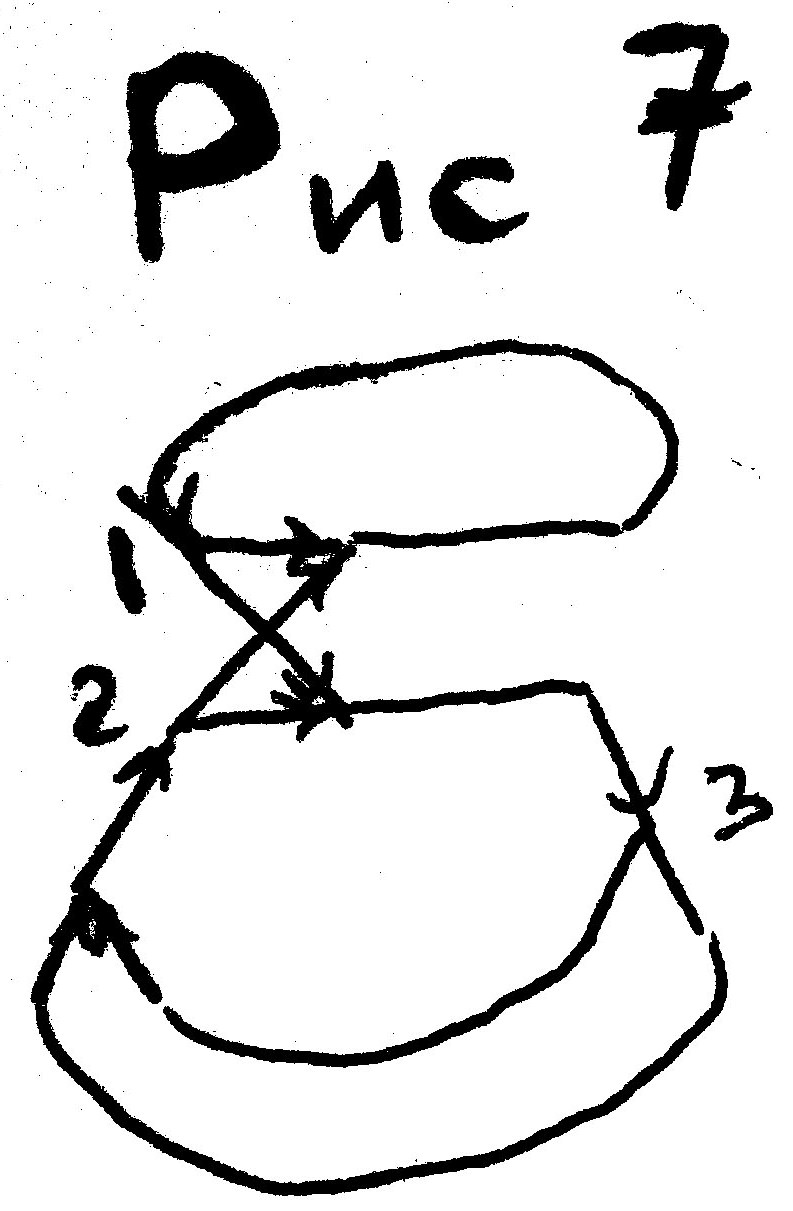}
	
	Легко понять, что при этом граф увеличится не более чем в 4 раза. Заглянем в таблицу \ref{table_second}, и осознаем, что нас интересует только случай когда за эти 4 операции произошло не больше одного запрета (в противном случае оценка стала верна). 
	
	\includegraphics[height=40mm]{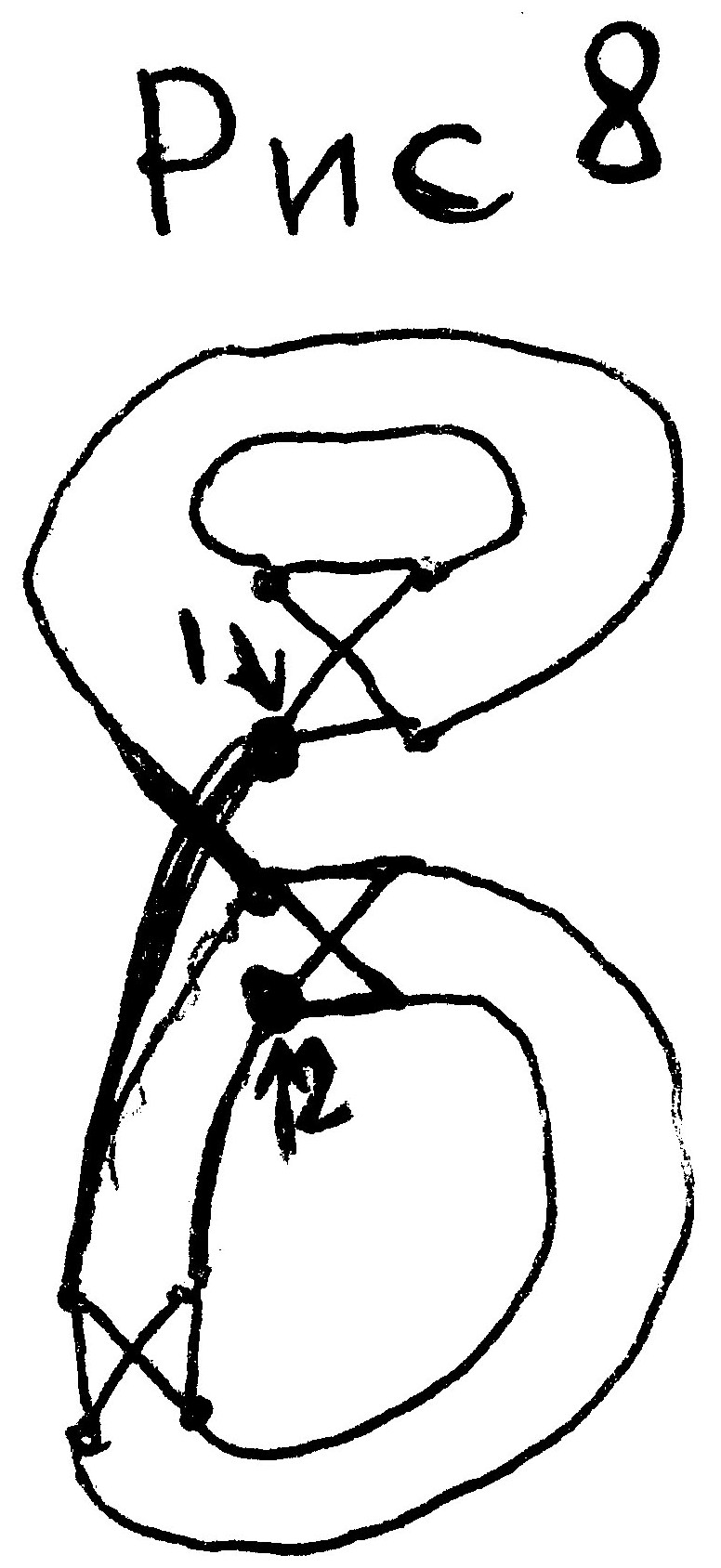}
	
	Теперь схлопнем те из развилок 1 и 2 (на последнем рисунке) которые можно. Это либо одна, либо две. Опять понятно что в итоге все равно будет не более чем учетверение (т.к. по каждому ребру исходного графа мы проехали в общей сложности не более чем 4 раза). Если это только одна, значит был запрет. 5 операций, $\leq + 4$ развилки, не более чем учетверение. См. таблицу \ref{table_second} --- оценка верна. Иначе --- 6 операций, не более чем учетверение --- оценка верна.
	
	Теперь второй случай. 
	
	Предположим, что развилки пронумерованы цифрами 1,2,3,4 в порядке слева направо как они располагаются на рисунке. 
	Схлопнем развилки 2 и 3. Запретов не было. Теперь схлопнем развилку 1 с её новой парой. Если запреты были, то 2 операции, не более чем +1 развилка, не более чем удвоение --- оценка верна. Если не было, то схлопнем еще и развилку 4 с её текущей парой, а затем те две развилки, одна из которых образовалась при этой операции, а вторая --- при схлопывании развилки 1. Если запреты были, то 4 операции, $\leq +3$ развилок, не более чем утроение --- оценка верна. Если запретов не было, то схлопнем теперь все имеющиеся в графе пары, которых три. Получим в итоге не более чем ушестерение, 7 операций, не более чем +7 развилок. Оценка верна.
	
\subsubsection{Неизменность количества развилок при схлапывании}

	Помним, что в этом пункте нам осталось только рассмотреть случай когда в графе только одна схлапываемая пара развилок. Если после первой операции схлапывания оставшаяся развилка (не важно, входящая или исходящая) --- попала на поддерево, где были еще развилки --- мы при первой операции двигаем развилку в соответствующую сторону и потому схлапываем новообразовавшуюся пару. Две операции, не более чем удвоение, запреты были --- оценка верна.
	
	Значит, остается случай когда с каждой стороны одно из ребер из корневой вершины ведет сразу в другую вершину, и на это ребро попадают обе развилки после первого схлапывания. 
	
	\includegraphics[height = 40mm]{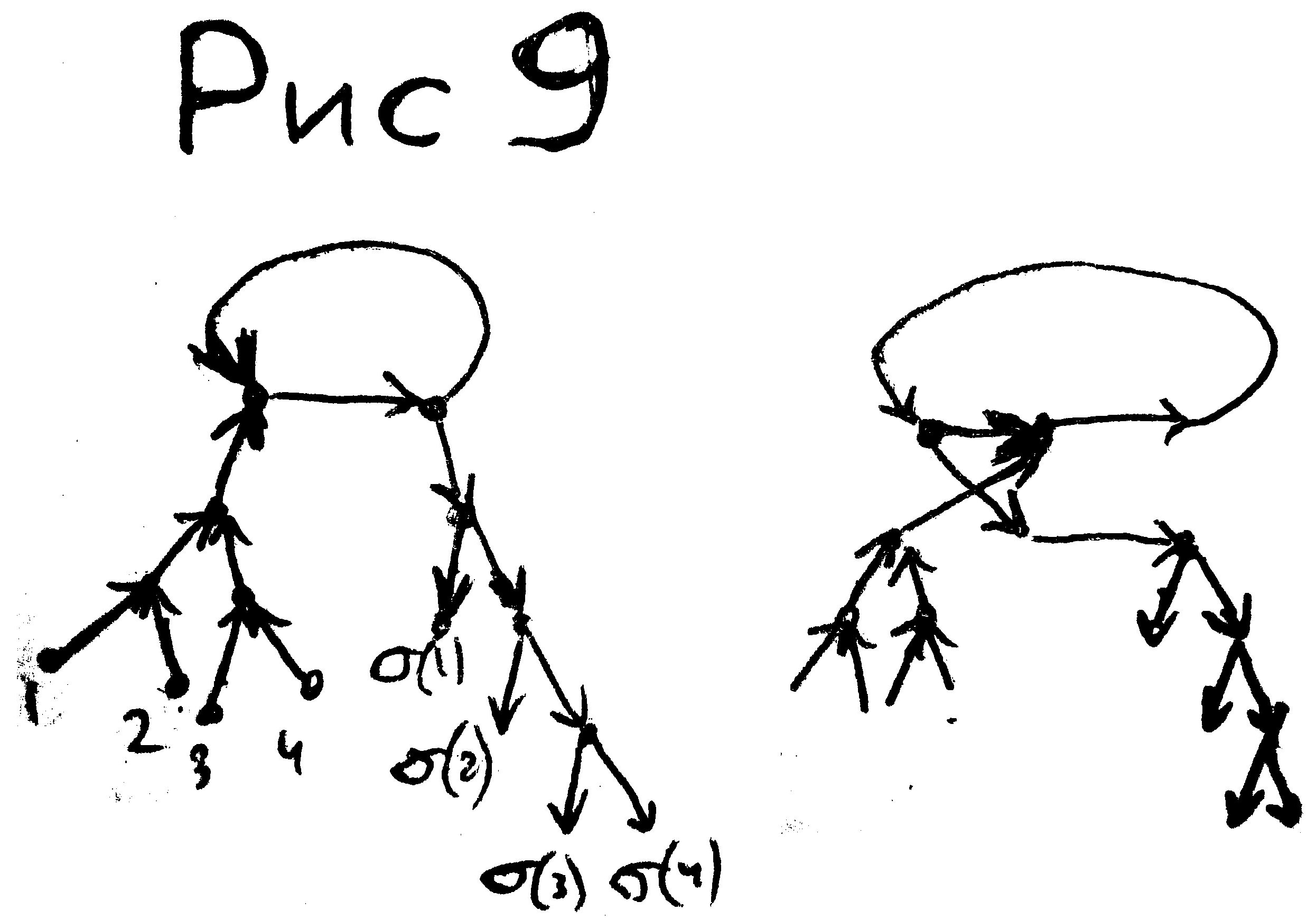} 
	
	Видно, что новый граф имеет структуру абсолютно аналогичную предыдущему. Ну снова схлопнем единственную имеющуюся развилку, и вообще будем продолжать это делать пока не случится иначе. А когда случится --- сделаем как описывалось ранее --- схлопнем с кем-нибудь из поддерева. Если это произошло $c$ раз, то граф увеличился не более чем в $c+1$ раз (надо понимать, что мы как бы заранее планировали с каким поддеревом будем схлапывать последнюю развилку, и в соответствующую сторону сдвигали нашу при схлопываниях). $c + 1$ операция, $c \geq 2$, не более чем $+ 2$ развилки. Оценка верна (очевидно из таблицы и анализа производной соотв. выражения) всегда кроме случая $c = 2$, $+ 2$ развилки. В этом случае схлопнем пару, аналогичную той которую уже схлопнули два раза (она есть, т.к. последих запретов не было --- отсюда было $+ 2$) --- сместились по таблице на 1 вправо вниз --- оценка верна. Тем более оценка верна, если в конце наших идентичных операций не +1 развилка, а -1 ($c$ операций, $-1$ развилка, увеличение не более чем в $c+1$ раз). Если вдруг после этого "`-1"' осталась только одна развилка --- то аналогично доказательствам приведенным ранее --- поскольку верна усиленная оценка на размер графа --- верно и утверждение нашей основной оценки для одноразвилковых графов.
	
\subsection{Начальное увеличение количества развилок}

	Пусть в $\mathbb{G}_k$ ровно одна развилка, и для него утверждение оценки верно. Если после схлопывания этой развилки не произойдет запретов, то количество развилок увеличится и станет равным двум, а нам будет необходимо сделать еще несколько операций так, чтобы асимптотика размера графа отстала в нужное количество раз или мы вернулись к одноразвилковой ситуации и утверждение оценки было верно. Без ограничения общности $\mathbb{G}_{k+1}$ будет выглядеть так:
	
	\includegraphics[height=20mm]{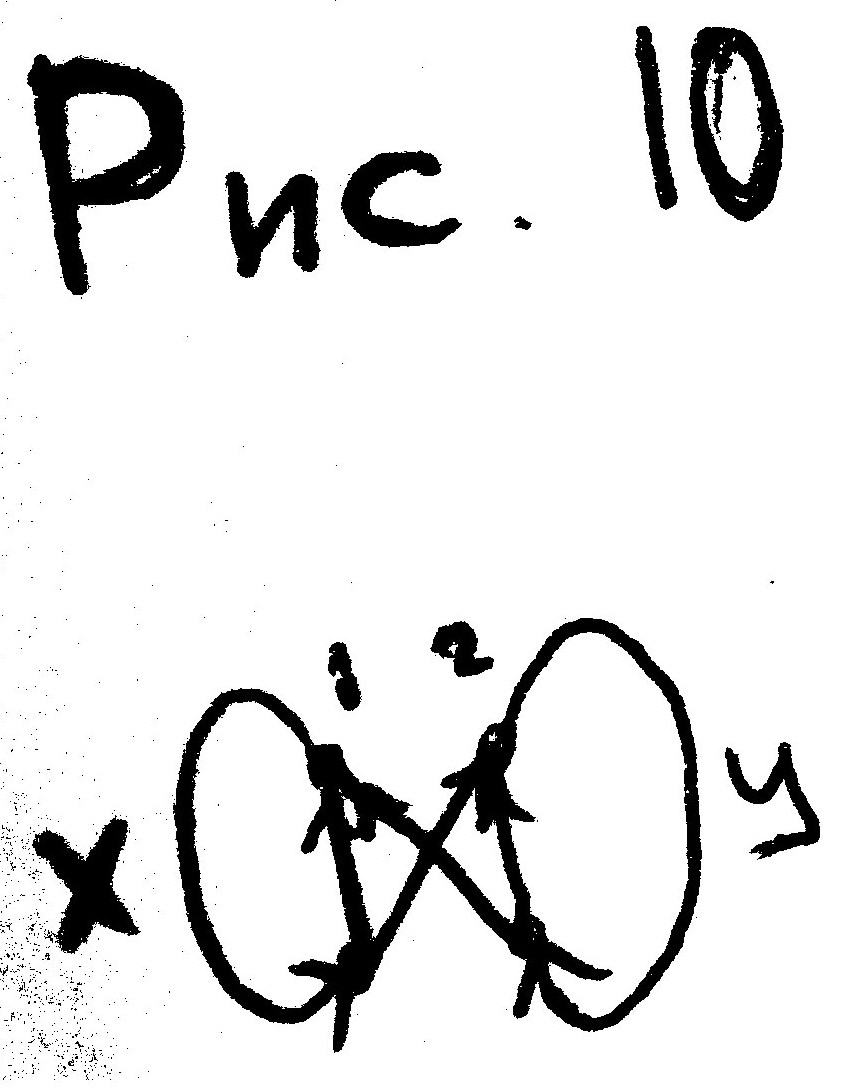}
	
	Нам известно, что $x + y \leq \varphi_{k+1}$, а $x + y + max(x,y) \leq \varphi_{k+2}$ (из утверждения оценки для $\mathbb G_k$). Схлопнем развилки 1 и 2. Получится такой граф (возмонжо, некоторые ребра из i,j,k и l запрещены --- но только они, иначе есть цикл со входящими развилками, но без входящих, или наоборот):
	
	\includegraphics[height=40mm]{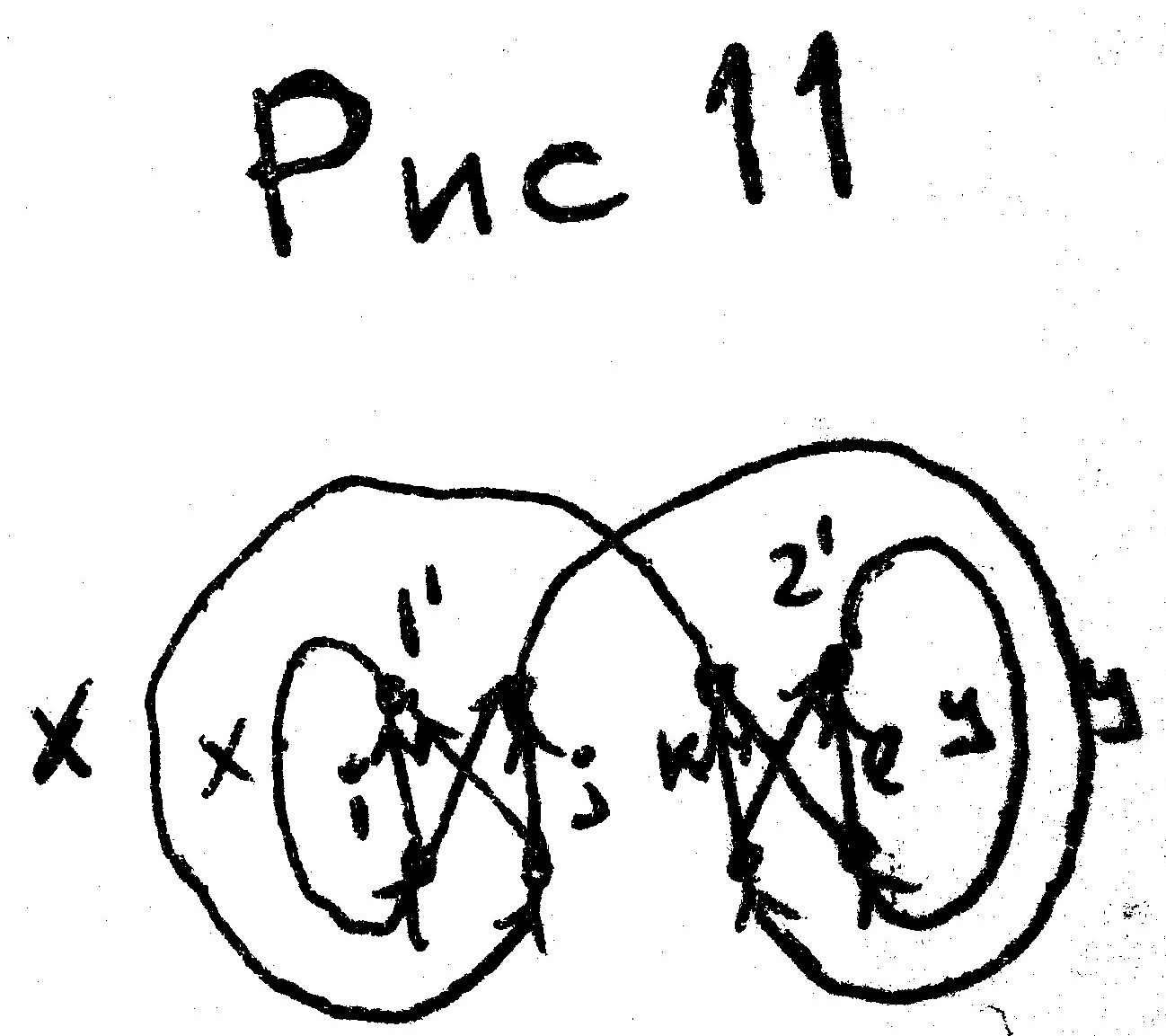}
	
	Если запрещены ребра i и j, то мы могли схлопнуть только развилку 1 и получить одноразвилковый граф с верным утверждением оценки. Аналогично если запрещены ребра k и l. Если осталась развилка $1'$ --- схлопнем её. Если $2'$ осталась --- её. Аналогично после этого появятся развилки 1'' или 2''. Будем продолжать схлопывать эти развилки пока они остаются. Когда-нибудь это закончится. Пусть с развилками $1^{(t)}$ это произведено $n_1$ раз, а с развилками $2^{(t)}$ $n_2$ раз. Граф увеличился не более чем в $max(n_1, n_2) + 1$ раз. Сделано $n_1 + n_2$ операций. Количество развилок --- стало не более чем $n_1 + n_2$ (вспомним, что оно определяет $k$ в таблице). Получается, что в таблице нас интересует второй столбец, когда $c = k$. Заметим, что $n_1, n_2 \geq 1$. Без ограничения общности можно предположить, что $n_1 \leq n_2$. Получается, что граф увеличился в $n_2 + 1$ раз. Тогда оценка не верна только для случаев $n_1 = 1, n_2 \in [1, 6]; n_1 = 2, n_2 \in [2, 3]$. Заметим, что если за эти $n_1 + n_2$ операций был хотя бы один запрет (кроме последних двух), то оценка не верна только когда $n_1 = 1, n_2 \in [2, 5]$, а если хотя бы два лишних запрета --- то в любом случае верна.
	
	Теперь схлопнем те из развилок 3' и 4', которые можно, а затем --- те из развилок 3'' и 4'' которые можно (на рисунке приведены результаты операций схлопывания развилок 1' и 2' --- $c$) и схлопывания 3' и 4' --- $d$) ). Теперь начнем оценивать числовые характеристики наших действий. 
	
	\includegraphics[height = 80mm]{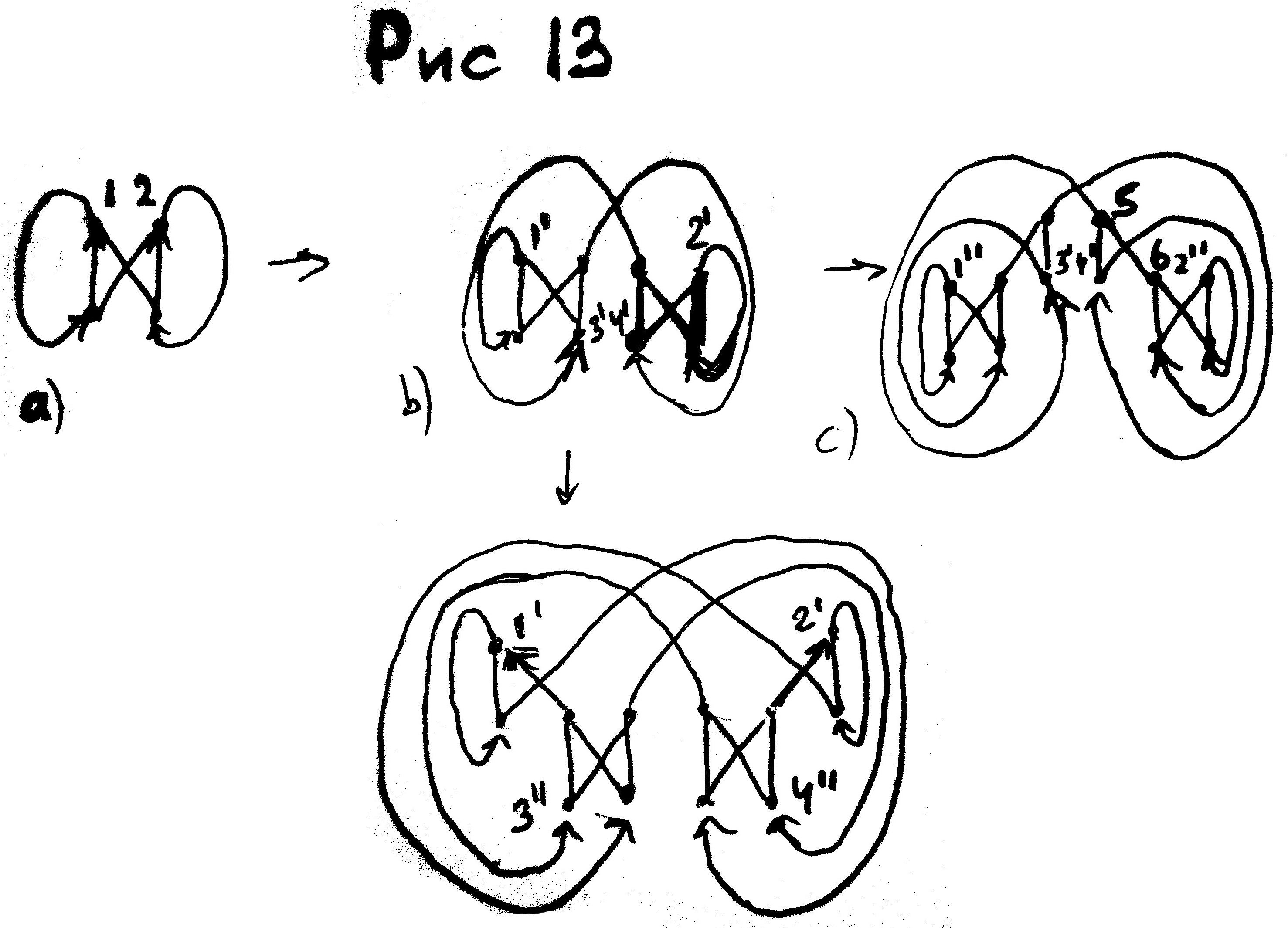}
	
	Если на первом этапе был ровно 1 запрет, то $n_1 = 1, n_2 \in [2, 5]$. Поскольку ребра i и j как мы условились одновременно быть запрещены не могли, то развилка 3' есть. Поскольку $n_2 \geq 2$, а запрет был ровно 1, то хотя бы одна из развилок 5, 6 есть. Значит можно с ней схлопнуть нашу развилку 3'. Числовые характеристики как если бы $n_1$ было равно 2 и был хотя бы один запрет, значит оценка верна.
	
	Если на первом этапе запретов не было --- значит все развилки есть. В частности --- 3' и 4'. Схлопнем их. Если запреты были, то числовые характеристики как если $n_1, n_2 \geq 2$, значит оценка верна. Если запретов не было опять, то есть 3'' и 4'', схлопнем их. Теперь числовые характеристики как если $n_1, n_2 \geq 3$, и оценка опять-таки верна. 
	
	Заметим, что если в итоге стала одна развилка, то для этого графа верна усиленная оценка, а если 0 --- то рассуждение окончено.

\subsection{Замечание}

	По сути данное доказательство --- есть формальное обоснование того факта, что если развилок много, то их встречи происходят очень часто относительно диаметра графа. Но поскольку в случае большого количества развилок структура графа становится непонятной и сложной, нельзя утверждать что не случится настолько удачной ситуации, что увеличение размера графа не окупит все предыдущие потери (а искуственно ситуацию когда на каждое размножение приходится по удвоению построить можно. Больше чем удвоение  --- нет, т.к. длина несамопересекающегося пути в графе очевидно не может превышать общего количества ребер в нем). Поэтому приходится делать такую усиленную оценку в зависимости от текущего количества перекрестков.
	
	Рассмотрим оставшийся случай, когда $\dfrac{\varphi_u}{\varphi_v * \alpha^{u-v}}$ все-таки равно минимальному значению 0.9270509831. Это возможно только когда мы на основании оценки для $|\mathbb G_2|$ пытаемся оценить $|\mathbb G_3|$. Если в $\mathbb{G}_3$ одна развилка --- утверждение верно (в $\mathbb{G}_3$ в любом случае 3 ребра), если нет --- мы, по нашему доказательству, делали дополнительные шаги (т.к. усиленная оценка была не верна --- т.к. верна только самая мягкая оценка --- для одноразвилкового случая).
	
\subsection{Заключение}

	Как мы доказали, как бы мы ни делали наши операции --- в конце получится циклический граф с количеством ребер равным периоду слова. Поскольку за одну операцию количество развилок в графе уменьшается не более чем на 1, а циклический граф --- граф без развилок, то на предпоследнем шаге в графе $\mathbb{G}_k$ была ровно одна развилка.  Значит для этого графа оценка была верна. Тогда $|\mathbb G_{k+1}| = 2x_k + y_k + z_k \leq \varphi_{k+1}$, в чем и заключается теорема \ref{mainth}.
	
\section{Доказательство минимальности оценки}

	Несложно строится пример последоваетльности слов с бесконечно возрастающим периодом задаваемых соответствующим минимально возможным количеством запретов: просто будем сохранять в графах Рози одну развилку и делать так, что $y_{k+1} = x_k, x_{k+1} = x_k + y_k$ --- тогда размер графа будет увеличиваться ровно по последовательности фибоначчи. На некотором шаге уничтожим обе развилки и получим соответствующее слово. Небольшие отклонения в этом процессе позволяют получить слова с другими периодами задаваемые соответствующим минимальным количеством запретов.
	
\section{Случай многобуквенного алфавита}	

	Случай многобуквенного алфавита довольно удачно сводится к случаю двухбуквенного. 
	Заметим, что в графах Рози слов содержащих k букв могут встречаться вершины входящей и исходящей степени от 1 до k. Аналогично двухбуквенному случаю все сводится к рассмотрению графов, в которых одна из входящей и исходящей степеней любой вершины равна 1. 
	
	Вершина входящей степени l может быть представлена как бинарное дерево (не важно какой формы) имеющее l концов (например, $l-1$ последовательных входящих развилок), каждое ребро которого имеет вес 0. 
	Легко понять как происходит эволюция графов Рози в случае k-буквенного алфавита: развилки аналогично едут навстречу друг другу, а при встрече, допустим, n-валентной входящей развилки и m-валентной исходящей, получается m n-валентных входящих развилок и n m-валентных исходящих. После этого как-то происходят запреты. 
	
	Заметим, что если сначала столкнуть развилки, а потом заменить на бинарные - получится $m * (n-1)$ входящих двоичных развилок, и $n * (m-1)$ исходящих, а если сначала заменить на деревья, а потом схлопнуть все эти развилки друг с другом, то получится столько же (каждая исходящая при прохождении через каждую входящюю дает +1 входящую развилку, тоесть всего $+ (m-1) * (n-1)$). Только вот запреты уже надо делать не в конце, а по ходу схлопывания двоичных - и тут их получится меньше, при той же итоговой конфигурации. 
	
	Получается, если в случае многобуквенного в начальном графе заменить многовалентные развилки на бинарные деревья, и провести эволюцию с аналогичными запрещениями до конца, то ребер будет столько же, а запретов меньше. Значит, максимальный конечный размер графа увеличивается не быстрее чем экспонента с основанием золотое сечение и степень - количество запретов, и интересна только начальная константа 
	(очевидно, из изложенного ранее рассуждения, что в начальном графе $k-1$ развилка, а его размер --- 1, таким образом начальная константа $\left(\frac{2}{\alpha}\right)^{k-1}$ подходит). 
	
	Пример же аналогично строится, если с самого начала убить все развилки (не совсем с самого - иначе это будет просто запрещение букв), кроме 1 входящей и 1 исходящей, и дальше сделать рост как последовательность Фибоначчи.
	
	Следовательно, возникает интересный вопрос: поскольку требование наличия k букв элементарно (по крайней мере на уровне идеи) сводится к случаю наличия 2 букв благодаря уничтожению большинства валентностей у развилок, интересно исследовать асимптотику роста максимального размера графа в зависимости от количества запретов при более строгих ограничениях, допустим, таких: суммарная валентность входящих развилок за вычетом их количества всегда не меньше $k$ (по сути - это количество эквивалентных двоичных развилок). 
	
	Она (асимптотика), очевидно, будет существенно меньше --- то есть для сложных слов требуется большое количество запретов. Впрочем, может быть, что она тоже будет экспоненциальной.

\section{Замечание.}
	Как мне стало известно от А.Я.Белова аналогичный результат был независимо получен И.И.Богдановым и Г.Р.Челноковым. Будет интересно сравнить доказательства когда появится их текст.

\appendix

\section{Таблица 1}

\label{table_pure}

	$\left(
	\begin{array}{cccccc} 
 \left(
\begin{array}{c}
 1.05902 \\
 2 \\
 2
\end{array}
\right) & \left(
\begin{array}{c}
 1.30902 \\
 2 \\
 1
\end{array}
\right) & \left(
\begin{array}{c}
 1.61803 \\
 2 \\
 0
\end{array}
\right) & \left(
\begin{array}{c}
 2. \\
 2 \\
 -1
\end{array}
\right) & \left(
\begin{array}{c}
 2.47214 \\
 2 \\
 -2
\end{array}
\right) & \left(
\begin{array}{c}
 3.05573 \\
 2 \\
 -3
\end{array}
\right) \\
 \left(
\begin{array}{c}
 1.38627 \\
 3 \\
 3
\end{array}
\right) & \left(
\begin{array}{c}
 1.71353 \\
 3 \\
 2
\end{array}
\right) & \left(
\begin{array}{c}
 2.11803 \\
 3 \\
 1
\end{array}
\right) & \left(
\begin{array}{c}
 2.61803 \\
 3 \\
 0
\end{array}
\right) & \left(
\begin{array}{c}
 3.23607 \\
 3 \\
 -1
\end{array}
\right) & \left(
\begin{array}{c}
 4. \\
 3 \\
 -2
\end{array}
\right) \\
 \left(
\begin{array}{c}
 1.81465 \\
 4 \\
 4
\end{array}
\right) & \left(
\begin{array}{c}
 2.24303 \\
 4 \\
 3
\end{array}
\right) & \left(
\begin{array}{c}
 2.77254 \\
 4 \\
 2
\end{array}
\right) & \left(
\begin{array}{c}
 3.42705 \\
 4 \\
 1
\end{array}
\right) & \left(
\begin{array}{c}
 4.23607 \\
 4 \\
 0
\end{array}
\right) & \left(
\begin{array}{c}
 5.23607 \\
 4 \\
 -1
\end{array}
\right) \\
 \left(
\begin{array}{c}
 2.37541 \\
 5 \\
 5
\end{array}
\right) & \left(
\begin{array}{c}
 2.93617 \\
 5 \\
 4
\end{array}
\right) & \left(
\begin{array}{c}
 3.62931 \\
 5 \\
 3
\end{array}
\right) & \left(
\begin{array}{c}
 4.48607 \\
 5 \\
 2
\end{array}
\right) & \left(
\begin{array}{c}
 5.54508 \\
 5 \\
 1
\end{array}
\right) & \left(
\begin{array}{c}
 6.8541 \\
 5 \\
 0
\end{array}
\right) \\
 \left(
\begin{array}{c}
 3.10945 \\
 6 \\
 6
\end{array}
\right) & \left(
\begin{array}{c}
 3.8435 \\
 6 \\
 5
\end{array}
\right) & \left(
\begin{array}{c}
 4.75082 \\
 6 \\
 4
\end{array}
\right) & \left(
\begin{array}{c}
 5.87234 \\
 6 \\
 3
\end{array}
\right) & \left(
\begin{array}{c}
 7.25861 \\
 6 \\
 2
\end{array}
\right) & \left(
\begin{array}{c}
 8.97214 \\
 6 \\
 1
\end{array}
\right) \\
 \left(
\begin{array}{c}
 4.07033 \\
 7 \\
 7
\end{array}
\right) & \left(
\begin{array}{c}
 5.0312 \\
 7 \\
 6
\end{array}
\right) & \left(
\begin{array}{c}
 6.21891 \\
 7 \\
 5
\end{array}
\right) & \left(
\begin{array}{c}
 7.68699 \\
 7 \\
 4
\end{array}
\right) & \left(
\begin{array}{c}
 9.50164 \\
 7 \\
 3
\end{array}
\right) & \left(
\begin{array}{c}
 11.7447 \\
 7 \\
 2
\end{array}
\right) \\
 \left(
\begin{array}{c}
 5.32813 \\
 8 \\
 8
\end{array}
\right) & \left(
\begin{array}{c}
 6.58593 \\
 8 \\
 7
\end{array}
\right) & \left(
\begin{array}{c}
 8.14065 \\
 8 \\
 6
\end{array}
\right) & \left(
\begin{array}{c}
 10.0624 \\
 8 \\
 5
\end{array}
\right) & \left(
\begin{array}{c}
 12.4378 \\
 8 \\
 4
\end{array}
\right) & \left(
\begin{array}{c}
 15.374 \\
 8 \\
 3
\end{array}
\right) \\
 \left(
\begin{array}{c}
 6.97461 \\
 9 \\
 9
\end{array}
\right) & \left(
\begin{array}{c}
 8.62109 \\
 9 \\
 8
\end{array}
\right) & \left(
\begin{array}{c}
 10.6563 \\
 9 \\
 7
\end{array}
\right) & \left(
\begin{array}{c}
 13.1719 \\
 9 \\
 6
\end{array}
\right) & \left(
\begin{array}{c}
 16.2813 \\
 9 \\
 5
\end{array}
\right) & \left(
\begin{array}{c}
 20.1248 \\
 9 \\
 4
\end{array}
\right) \\
 \left(
\begin{array}{c}
 9.12988 \\
 10 \\
 10
\end{array}
\right) & \left(
\begin{array}{c}
 11.2852 \\
 10 \\
 9
\end{array}
\right) & \left(
\begin{array}{c}
 13.9492 \\
 10 \\
 8
\end{array}
\right) & \left(
\begin{array}{c}
 17.2422 \\
 10 \\
 7
\end{array}
\right) & \left(
\begin{array}{c}
 21.3125 \\
 10 \\
 6
\end{array}
\right) & \left(
\begin{array}{c}
 26.3437 \\
 10 \\
 5
\end{array}
\right) \\
 \left(
\begin{array}{c}
 11.9512 \\
 11 \\
 11
\end{array}
\right) & \left(
\begin{array}{c}
 14.7725 \\
 11 \\
 10
\end{array}
\right) & \left(
\begin{array}{c}
 18.2598 \\
 11 \\
 9
\end{array}
\right) & \left(
\begin{array}{c}
 22.5703 \\
 11 \\
 8
\end{array}
\right) & \left(
\begin{array}{c}
 27.8984 \\
 11 \\
 7
\end{array}
\right) & \left(
\begin{array}{c}
 34.4844 \\
 11 \\
 6
\end{array}
\right)
\end{array}
	\right)$

\section{Таблица 2}
\label{table_evasion}

	$\left(
\begin{array}{cccccccc}
 1.000000 & 0.8090170 & 0.8726780 & 0.8472136 & 0.8567627 & 0.8530900 & 0.8544891 & 0.8539542 \\
 1.236068 & 1.000000 & 1.078689 & 1.047214 & 1.059017 & 1.054477 & 1.056207 & 1.055545 \\
 1.145898 & 0.9270510 & 1.000000 & 0.9708204 & 0.9817627 & 0.9775541 & 0.9791574 & 0.9785444 \\
 1.180340 & 0.9549150 & 1.030057 & 1.000000 & 1.011271 & 1.006936 & 1.008588 & 1.007956 \\
 1.167184 & 0.9442719 & 1.018576 & 0.9888544 & 1.000000 & 0.9957132 & 0.9973463 & 0.9967219 \\
 1.172209 & 0.9483372 & 1.022961 & 0.9931116 & 1.004305 & 1.000000 & 1.001640 & 1.001013 \\
 1.170290 & 0.9467844 & 1.021286 & 0.9914855 & 1.002661 & 0.9983626 & 1.000000 & 0.9993739 \\
 1.171023 & 0.9473775 & 1.021926 & 0.9921066 & 1.003289 & 0.9989880 & 1.000626 & 1.000000
\end{array}
\right)$

\section{Таблица 3}
\label{table_second}

	$\left(
	\begin{array}{cccccc}
 \left(
\begin{array}{c}
 1. \\
 2 \\
 2
\end{array}
\right) & \left(
\begin{array}{c}
 1.23607 \\
 2 \\
 1
\end{array}
\right) & \left(
\begin{array}{c}
 1.52786 \\
 2 \\
 0
\end{array}
\right) & \left(
\begin{array}{c}
 1.88854 \\
 2 \\
 -1
\end{array}
\right) & \left(
\begin{array}{c}
 2.33437 \\
 2 \\
 -2
\end{array}
\right) & \left(
\begin{array}{c}
 2.88544 \\
 2 \\
 -3
\end{array}
\right) \\
 \left(
\begin{array}{c}
 1.30902 \\
 3 \\
 3
\end{array}
\right) & \left(
\begin{array}{c}
 1.61803 \\
 3 \\
 2
\end{array}
\right) & \left(
\begin{array}{c}
 2. \\
 3 \\
 1
\end{array}
\right) & \left(
\begin{array}{c}
 2.47214 \\
 3 \\
 0
\end{array}
\right) & \left(
\begin{array}{c}
 3.05573 \\
 3 \\
 -1
\end{array}
\right) & \left(
\begin{array}{c}
 3.77709 \\
 3 \\
 -2
\end{array}
\right) \\
 \left(
\begin{array}{c}
 1.71353 \\
 4 \\
 4
\end{array}
\right) & \left(
\begin{array}{c}
 2.11803 \\
 4 \\
 3
\end{array}
\right) & \left(
\begin{array}{c}
 2.61803 \\
 4 \\
 2
\end{array}
\right) & \left(
\begin{array}{c}
 3.23607 \\
 4 \\
 1
\end{array}
\right) & \left(
\begin{array}{c}
 4. \\
 4 \\
 0
\end{array}
\right) & \left(
\begin{array}{c}
 4.94427 \\
 4 \\
 -1
\end{array}
\right) \\
 \left(
\begin{array}{c}
 2.24303 \\
 5 \\
 5
\end{array}
\right) & \left(
\begin{array}{c}
 2.77254 \\
 5 \\
 4
\end{array}
\right) & \left(
\begin{array}{c}
 3.42705 \\
 5 \\
 3
\end{array}
\right) & \left(
\begin{array}{c}
 4.23607 \\
 5 \\
 2
\end{array}
\right) & \left(
\begin{array}{c}
 5.23607 \\
 5 \\
 1
\end{array}
\right) & \left(
\begin{array}{c}
 6.47214 \\
 5 \\
 0
\end{array}
\right) \\
 \left(
\begin{array}{c}
 2.93617 \\
 6 \\
 6
\end{array}
\right) & \left(
\begin{array}{c}
 3.62931 \\
 6 \\
 5
\end{array}
\right) & \left(
\begin{array}{c}
 4.48607 \\
 6 \\
 4
\end{array}
\right) & \left(
\begin{array}{c}
 5.54508 \\
 6 \\
 3
\end{array}
\right) & \left(
\begin{array}{c}
 6.8541 \\
 6 \\
 2
\end{array}
\right) & \left(
\begin{array}{c}
 8.47214 \\
 6 \\
 1
\end{array}
\right) \\
 \left(
\begin{array}{c}
 3.8435 \\
 7 \\
 7
\end{array}
\right) & \left(
\begin{array}{c}
 4.75082 \\
 7 \\
 6
\end{array}
\right) & \left(
\begin{array}{c}
 5.87234 \\
 7 \\
 5
\end{array}
\right) & \left(
\begin{array}{c}
 7.25861 \\
 7 \\
 4
\end{array}
\right) & \left(
\begin{array}{c}
 8.97214 \\
 7 \\
 3
\end{array}
\right) & \left(
\begin{array}{c}
 11.0902 \\
 7 \\
 2
\end{array}
\right) \\
 \left(
\begin{array}{c}
 5.0312 \\
 8 \\
 8
\end{array}
\right) & \left(
\begin{array}{c}
 6.21891 \\
 8 \\
 7
\end{array}
\right) & \left(
\begin{array}{c}
 7.68699 \\
 8 \\
 6
\end{array}
\right) & \left(
\begin{array}{c}
 9.50164 \\
 8 \\
 5
\end{array}
\right) & \left(
\begin{array}{c}
 11.7447 \\
 8 \\
 4
\end{array}
\right) & \left(
\begin{array}{c}
 14.5172 \\
 8 \\
 3
\end{array}
\right) \\
 \left(
\begin{array}{c}
 6.58593 \\
 9 \\
 9
\end{array}
\right) & \left(
\begin{array}{c}
 8.14065 \\
 9 \\
 8
\end{array}
\right) & \left(
\begin{array}{c}
 10.0624 \\
 9 \\
 7
\end{array}
\right) & \left(
\begin{array}{c}
 12.4378 \\
 9 \\
 6
\end{array}
\right) & \left(
\begin{array}{c}
 15.374 \\
 9 \\
 5
\end{array}
\right) & \left(
\begin{array}{c}
 19.0033 \\
 9 \\
 4
\end{array}
\right) \\
 \left(
\begin{array}{c}
 8.62109 \\
 10 \\
 10
\end{array}
\right) & \left(
\begin{array}{c}
 10.6563 \\
 10 \\
 9
\end{array}
\right) & \left(
\begin{array}{c}
 13.1719 \\
 10 \\
 8
\end{array}
\right) & \left(
\begin{array}{c}
 16.2813 \\
 10 \\
 7
\end{array}
\right) & \left(
\begin{array}{c}
 20.1248 \\
 10 \\
 6
\end{array}
\right) & \left(
\begin{array}{c}
 24.8756 \\
 10 \\
 5
\end{array}
\right) \\
 \left(
\begin{array}{c}
 11.2852 \\
 11 \\
 11
\end{array}
\right) & \left(
\begin{array}{c}
 13.9492 \\
 11 \\
 10
\end{array}
\right) & \left(
\begin{array}{c}
 17.2422 \\
 11 \\
 9
\end{array}
\right) & \left(
\begin{array}{c}
 21.3125 \\
 11 \\
 8
\end{array}
\right) & \left(
\begin{array}{c}
 26.3437 \\
 11 \\
 7
\end{array}
\right) & \left(
\begin{array}{c}
 32.5626 \\
 11 \\
 6
\end{array}
\right)
\end{array}
	\right)$

\end{document}